\documentclass[final,1p,times]{elsarticle}

\usepackage{amssymb}
\usepackage{amsmath} \usepackage{amsthm}
\usepackage{lineno}
\usepackage{xcolor}
\usepackage{geometry}
\newgeometry{left=1.25in, right=1.25in}

\usepackage[utf8]{inputenc}
\usepackage[english]{babel}
\usepackage{tikz}
\usetikzlibrary{shapes.geometric}
\usetikzlibrary{arrows.meta,positioning,fit,calc}

\usepackage{graphicx}
\usepackage{caption, subcaption, float}
\usepackage[colorlinks=true, allcolors=blue]{hyperref}
\usepackage{verbatim}
\usepackage{algpseudocode}
\usepackage[linesnumbered,ruled,vlined]{algorithm2e} 
\usepackage{soul}
\usepackage{enumitem}
\usepackage{amsmath} 
\usepackage{cleveref}

\theoremstyle{plain}

\theoremstyle{definition}

\theoremstyle{remark}
\newtheorem*{remark*}{Remark}        
\usepackage[colorlinks=true, allcolors=blue]{hyperref}
\usepackage{algorithmicx}
\usepackage{algcompatible}
\usepackage{multirow}
\usepackage{booktabs}

\newlist{computesteps}{description}{1}
\setlist[computesteps]{labelsep=0.75em,leftmargin=*,font=\ttfamily}

\usepackage{tikz}
\usepackage{amsmath,amssymb,bm}
\usetikzlibrary{
  arrows.meta,
  backgrounds,
  calc,
  fit,
  positioning,
  shadows.blur
}

\usepackage[acronym]{glossaries}
\setkeys{glslink}{hyper=false}        
\newacronym{pod}{POD}{Proper Orthogonal Decomposition}
\newacronym{nn}{NN}{Neural Network}

\journal{journal name}

\begin{document}
\begin{frontmatter}

\title{SparseModesNet for coupled thermal-hydraulic-mechanical system \\ in heterogeneous porous media} 


\author[fsu]{Kapil Chawla} 
\author[gt]{Tomoki Koike}
\author[fsu]{Sanghyun Lee} 
\author[nl]{Julie Bessac}

\affiliation[fsu]{organization={Department of Mathematics, Florida State University},
            addressline={1017 Academic Way}, 
            city={Tallahassee},
            postcode={32306}, 
            state={FL},
            country={USA}}

\affiliation[gt]{organization={School of Aerospace Engineering, Georgia Institute of Technology},
            addressline={270 Ferst Drive NW}, 
            city={Atlanta},
            postcode={30332}, 
            state={GA},
            country={USA}}

\affiliation[nl]{organization={Computational Science Center, National Laboratory of the Rockies (NLR)},
            addressline={15013 Denver West Parkway}, 
            city={Golden},
            postcode={80401}, 
            state={CO},
            country={USA}}

\begin{abstract} 
The coupled thermo-hydro-mechanical (THM) models describe the interaction
between displacement, fluid pressure, and temperature in geothermal systems. Solving this coupled system on a fine mesh over many time steps is
computationally expensive because a large number of unknowns must be computed at each time step. The problem becomes more challenging when some of the coefficients such as permeability and thermal conductivity are highly
heterogeneous and when fractures create complex paths for fluid flow and heat
transfer. Reduced-order models (ROMs) approximate a full solution containing many
unknowns using only a small number of variables. One common method for
constructing a ROM is Proper Orthogonal Decomposition (POD), which finds the
dominant patterns in the solution, called modes. However, because POD uses a
linear representation, it may require many modes to accurately reconstruct
complex THM solutions. In this work, we apply SparseModesNet, a
method developed in~\cite{koike2026SparseModesNet} that parametrizes a
nonlinear decoder by a neural network (NN), to obtain higher reconstruction
accuracies of the coupled THM solution. During the NN training, SparseModesNet
starts from a larger set of candidate modes and selects a small number of
modes that are optimal for reducing the reconstruction error but are not
necessarily highest energy modes. The trained SparseModesNet then reconstructs the
full displacement, pressure, and temperature fields from the selected reduced
coordinates. We compare SparseModesNet with standard POD using three numerical
experiments. These experiments include a two-dimensional problem with
heterogeneous permeability and thermal conductivity, a two-dimensional connected
fracture network, and a three-dimensional problem based on a Utah FORGE fracture
geometry. We also compare pressure and temperature quantities of interest
obtained from the reconstructed fields. In addition, we test temporal
interpolation by training SparseModesNet on a subset of the solution snapshots
and reconstructing the remaining snapshots within the same time interval. The
results show that SparseModesNet produces much smaller reconstruction error than
POD for the same reduced dimension. It also accurately reproduces the coupled
solution fields and the selected geothermal quantity of interest.
\end{abstract}



\end{frontmatter}

\section{Introduction}
\label{sec:introduction}

Enhanced geothermal systems (EGSs) produce energy from hot rocks deep below the Earth's surface. Cold fluid is injected into the reservoir, where it flows through natural or artificial fractures, absorbs heat from the surrounding rock, and returns to the surface through production wells. During this circulation, changes in fluid pressure and temperature can also deform the surrounding rock. To account for these interactions, geothermal reservoirs are commonly modeled using coupled thermo-hydro-mechanical (THM) models~\cite{biot1941consolidation,mctigue1986thermoelastic,coussy2004poromechanics}. Solving a coupled THM model on a fine spatial mesh is challenging since it involves a large number of degrees of freedom at each time step, even more so when permeability and thermal conductivity vary throughout the reservoir. This heterogeneity is what creates complex spatial variations in the coupled THM solution~\cite{rutqvist2003hydromechanical,tsang2024coupled,podgorney2021forge}.

Many downstream geothermal tasks, including uncertainty quantification, design optimization, and long-horizon forecasting, require this expensive solve to be repeated many times, which motivates the construction of cheaper surrogates through, for example, projection-based reduced-order modeling (ROM)~\cite{benner2015survey,quarteroni2016reduced,hesthaven2016certified,benner2017model}. Projection-based ROM proceeds in two stages. First, a low-dimensional representation of the solution is identified from a set of solution snapshots, where the THM solution saved at each time step is called a solution snapshot. Second, a reduced dynamical model is constructed to evolve the coordinates of that representation. The accuracy attainable in the second stage is bounded by the first, because a representation that cannot reconstruct the snapshots at low dimension cannot support an accurate reduced model. For this reason, the reduced representation is the stage that determines whether an efficient ROM is possible. In this work we focus on this first stage, that is, on finding a compact and accurate representation of coupled THM solutions; the construction of the reduced dynamics is a separate step that we leave to future work.

By far the most widely-used method for state compression in model reduction literature is the proper orthogonal decomposition (POD) representation~\cite{benner2015survey,kramer2026reduced}. POD is a linear representation, in which the system state is approximated as a linear combination of a small number of POD modes. These modes are the leading left singular vectors of a data matrix consisting of state snapshots. The term POD originated in the analysis of turbulent flows~\cite{sirovich1987turbulence, berkooz1993proper,rowley2017model}. While, in other fields, the POD modes are sometimes called the principal components of the data or the Karhunen-Lo\'eve modes~\cite{breuer1991use,winter1992eigenfunction,romanowski1996reduced}.  The POD state representation can be used to derive a compact reduced dynamical model in several ways, including intrusive (requiring access to the code implementing the high-dimensional model operators) Petrov-Galerkin projection~\cite{iollo2000stability,lorenzi2016pod,ullmann2016pod}, as well as popular non-intrusive data-driven methods such as those based on the dynamic mode decomposition (DMD)~\cite{schmid2010dynamic,rowley2009spectral,tu2014dynamic} or on the Operator Inference (OpInf) framework~\cite{peherstorfer2016data,kramer2024learning}.

POD methods are popular because of their general applicability as well as due to the linear optimality property of POD: the POD modes span the optimal linear subspace in which the reconstruction error of the snapshot data is minimized~\cite{eckart1936approximation}. 
With this error guarantee, POD has proven effective in geothermal applications, for example, in the reduced-order modeling of linear thermo-poroelasticity~\cite{BallarinLeeYi2024,florez2020global} and fractured geothermal reservoirs~\cite{li2020geothermalpod}. However, even optimal linear approximations face challenges in efficient representation of states for systems exhibiting complex nonlinear dynamics; in the model reduction literature this phenomenon is known as the challenge of slow decay of the Kolmogorov $n$-width~\cite{maday2002Global,cohen2016Kolmogorov,peherstorfer2022Breaking}. Specifically, in the THM problem considered here, the medium is strongly heterogeneous: permeability and thermal conductivity change sharply between the fractures and the surrounding matrix, and this heterogeneity ingrains fine, spatially varying structure on the displacement, pressure, and temperature fields. The heterogeneous medium also channels flow through the factures, producing thermal and pressure fronts that propagate along the fracture network. Both effects, a spatially heterogeneous solution structure and transport-dominated fronts, are associated with slow Kolmogorov \(n\)-width decay, so a linear representation may require many leading modes to reach a low reconstruction error. To overcome this challenge, several \textit{nonlinear} methods for compressing the system state have been explored, including polynomial POD manifolds that use polynomial functions of the POD coefficients to reconstruct the state~\cite{geelen2023operator,geelen2023Learning,barnett2022quadratic,schwerdtner2024greedy}, fully deep-learning approaches that use neural networks (NNs) to both encode the state in a compact representation and decode the state again~\cite{lee2020model,kadeethum2021poroelasticity,kadeethum2022reduced,ballini2024application}, and hybrid approaches that use POD to encode the state and a neural network to decode the state~\cite{koike2026SparseModesNet,barnett2023Neural,aresdeparga2026Nonlinear}. 

In this work, we use the hybrid approach proposed in~\cite{koike2026SparseModesNet} called SparseModesNet (SMN). SMN begins with $s$ candidate POD modes and selects $r < s$ modes that lead to the lowest reconstruction error, where the selected modes are not required to be the leading POD modes. The coefficients of the selected modes are passed to a nonlinear decoder parameterized by a NN, which reconstructs the fine-mesh THM solution. The mode selection addresses the second limitation and the nonlinear decoder addresses the first, so that SMN can reach a given reconstruction accuracy at a smaller reduced dimension than POD.

In this work, we apply SMN to obtain a compact nonlinear representation of coupled THM solutions, involving displacement, pressure, and temperature, in domains with heterogeneous permeability and thermal conductivity and complex fracture networks. This is the heterogeneity-dominated regime identified above: the sharp, fracture-controlled variation that makes a global linear basis inefficient also stresses a learned nonlinear representation, so accurate reduction in this setting cannot be assumed. To our knowledge, nonlinear manifold reduction has seen limited demonstration on coupled THM solutions in reservoirs with this degree of heterogeneity and fracture complexity, and the quantities that matter for geothermal engineering are specific outputs rather than the global field error alone. A central contribution of this work is therefore to show that SMN yields an accurate low-dimensional representation of these strongly heterogeneous THM dynamics, and preserves the engineering quantities of interest, whereas standard POD at the same number of modes fails to achieve the same order of accuracy.
 
Our numerical experiments include a randomly heterogeneous coefficient field, a two-dimensional fracture network, and a three-dimensional fracture geometry based on the Utah FORGE site~\cite{finnila2025utahforge_dfn}. 
The Utah FORGE is a DOE-funded underground field laboratory in Utah, dedicated to advancing EGS by testing and refining technologies to extract heat from deep, low-permeability rock formations. Its goal is to create a replicable, commercial pathway for EGS, making geothermal energy viable in new regions worldwide. 
We measure the reconstruction error for the entire THM solution snapshots and separately in displacement, pressure, and temperature. We also examine whether the reconstructed fields preserve important geothermal quantities of interest, including the average temperature near the production wells and the injection-production pressure difference normalized by the circulation rate. Finally, we evaluate temporal interpolation by training SMN on a subset of the THM snapshots and testing it on the remaining unseen snapshots, which probes whether the representation generalizes across the trajectory in terms of interpolation at minimum, a prerequisite for its later use in a predictive reduced model.

The remainder of this paper is organized as follows. Section~\ref{sec:thm_model} introduces the coupled THM model and its numerical discretization. Section~\ref{sec:SparseModesNet} describes POD, SparseModesNet, and the reconstruction-error measures. Section~\ref{sec:numerical} presents the three numerical experiments, the corresponding quantities of interest, and the temporal interpolation results. Finally, Section~\ref{sec:conclusion} summarizes the main findings and discusses possible extensions.

\section{Thermo-Hydro-Mechanical Model Formulation}
\label{sec:thm_model}

In this section, we consider a linear thermo-poroelasticity model to describe coupled thermo-hydro-mechanical (THM) processes in a deformable porous medium. The model couples the mechanical displacement of the solid, the fluid pressure, and the temperature field through momentum balance, fluid mass conservation, and energy balance equations. This formulation follows the standard linearized thermo-poroelasticity framework considered in \cite{BallarinLeeYi2024}. The finite element solutions of this model at
different time levels are referred to as the high-fidelity solution snapshots
and serve as the reference data for the SparseModesNet approximation.

\begin{figure}[!h]
    \centering
    \begin{subfigure}[t]{0.48\textwidth}
        \centering
        \includegraphics[width=\textwidth]{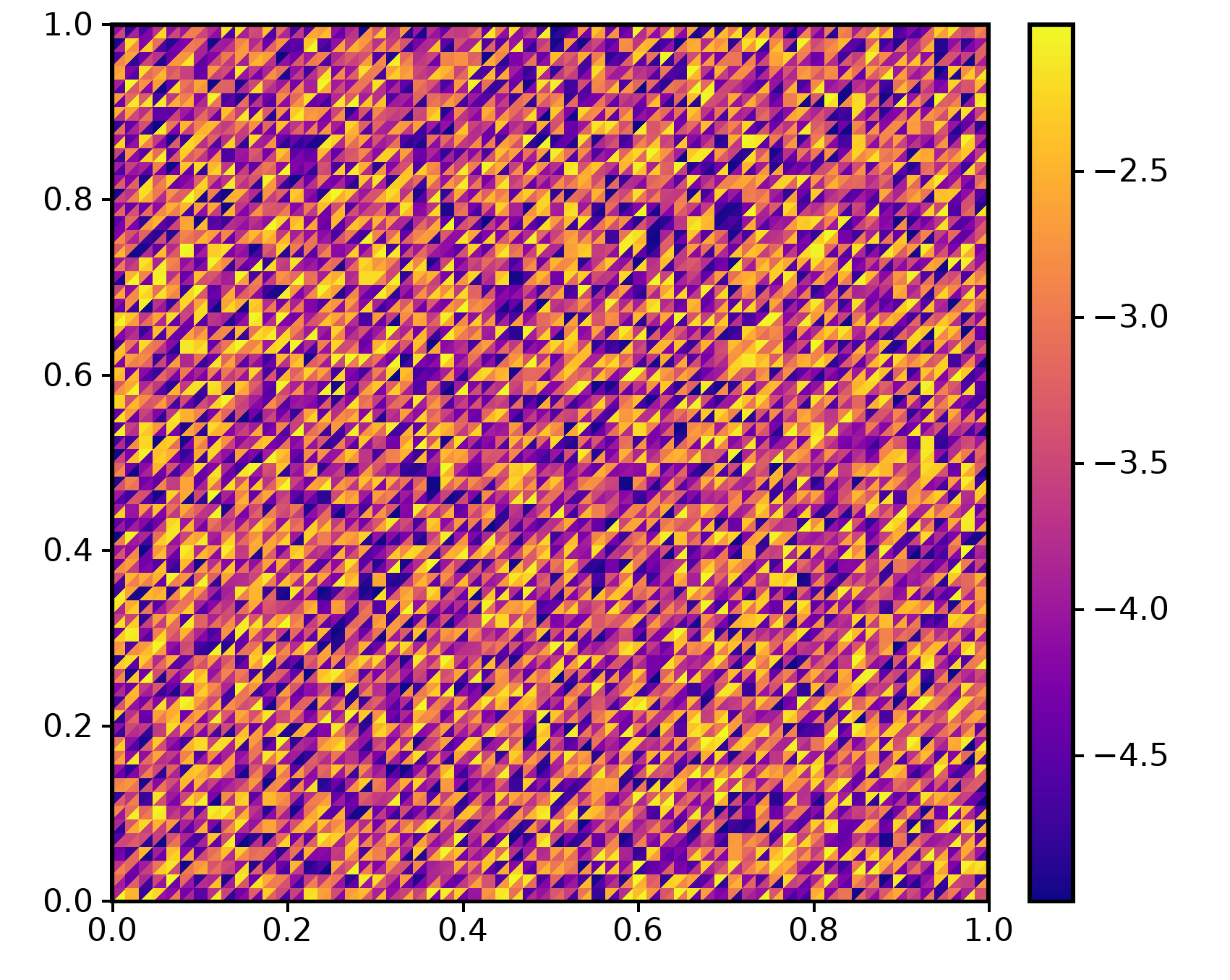}
        \caption{Permeability coefficient $\kappa(\mathbf{x})$}
        \label{fig:log10_kappa_field}
    \end{subfigure}
    \hfill
    \begin{subfigure}[t]{0.45\textwidth}
        \centering
        \includegraphics[width=\textwidth]{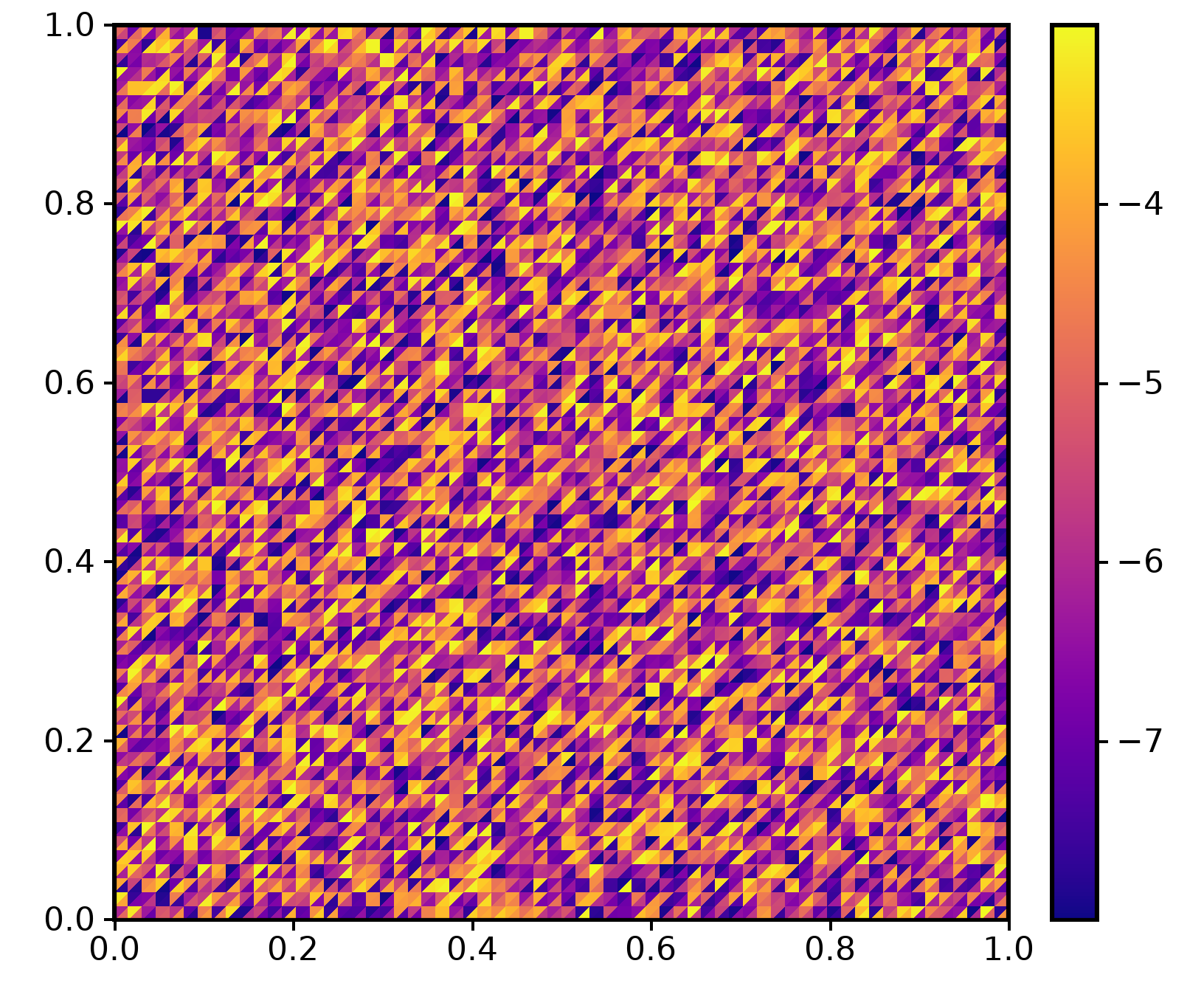}
        \caption{Thermal conductivity coefficient $d_T(\mathbf{x})$}
        \label{fig:log10_conductivity_field}
    \end{subfigure}
    \caption{Representative heterogeneous coefficient fields used in the THM numerical experiments. The permeability coefficient $\kappa(\mathbf{x})$ and thermal conductivity coefficient $d_T(\mathbf{x})$ are displayed on a logarithmic scale.}
    \label{fig:random_KD_fields}
\end{figure}

Let $\Omega \subset \mathbb{R}^{d}$, $d=2,3$, be a bounded, connected, Lipschitz domain with boundary $\partial\Omega$, and let $\mathcal{J}=(0,T]$ denote the time interval. We seek the displacement $\mathbf{u}:\Omega\times \mathcal{J}\rightarrow \mathbb{R}^{d}$, the fluid pressure $p:\Omega\times \mathcal{J}\rightarrow \mathbb{R}$, and the temperature field $\theta:\Omega\times \mathcal{J}\rightarrow \mathbb{R}$ satisfying

\begin{subequations}
\label{eq:thm_strong}
\begin{align}
-\nabla \cdot
\left(
\boldsymbol{\sigma}(\mathbf{u})
-
\alpha p \mathbf{I}
-
3\alpha_T K_{\mathrm{dr}} \theta \mathbf{I}
\right)
&=
\mathbf{f}
&& \text{in } \Omega \times \mathcal{J},
\label{eq:thm_momentum}
\\
\frac{\partial}{\partial t}
\left(
c_0 p
+
\alpha \nabla \cdot \mathbf{u}
-
3\alpha_m \theta
\right)
-
\nabla \cdot
\left(
\mathbf{K}\nabla p
\right)
&=
g
&& \text{in } \Omega \times \mathcal{J},
\label{eq:thm_pressure}
\\
\frac{\partial}{\partial t}
\left(
C_d \theta
+
3\alpha_T K_{\mathrm{dr}}\theta_0 \nabla \cdot \mathbf{u}
-
3\alpha_m \theta_0 p
\right)
-
\nabla \cdot
\left(
\mathbf{D}\nabla \theta
\right)
&=
\eta
&& \text{in } \Omega \times \mathcal{J}.
\label{eq:thm_temperature}
\end{align}
\end{subequations}
Here, the linear elastic stress tensor is defined by
$\boldsymbol{\sigma}(\mathbf{u})
=
2\mu \boldsymbol{\varepsilon}(\mathbf{u})
+
\lambda(\nabla\cdot \mathbf{u})\mathbf{I}_d$,
where $\lambda$ and $\mu$ are the Lamé parameters, and
$\boldsymbol{\varepsilon}(\mathbf{u})
=
\frac{1}{2}
\left(
\nabla \mathbf{u}
+
\nabla \mathbf{u}^{T}
\right)$
is the infinitesimal strain tensor. The drained bulk modulus is given by
$K_{\mathrm{dr}}=(d\lambda+2\mu)/d$. The terms $\mathbf{f}$, $g$, and $\eta$ denote the body force, the volumetric fluid source/sink term, and the volumetric heat source/sink term, respectively. The parameter $\alpha$ is the Biot--Willis coefficient, $c_0$ is the constrained specific storage coefficient, $\alpha_T$ is the volumetric thermal dilation coefficient, $\alpha_m$ is the thermal dilation coefficient, $C_d$ is the effective heat capacity, and $\theta_0$ is a reference temperature. Finally, in this work we consider isotropic heterogeneous permeability and thermal conductivity fields. For $\mathbf{x}\in\Omega$, these fields are written as
$$
\mathbf{K}(\mathbf{x})=\kappa(\mathbf{x})\mathbf{I}_d,
\qquad
\mathbf{D}(\mathbf{x})=d_T(\mathbf{x})\mathbf{I}_d,
$$
where $\kappa$ and $d_T$ are scalar coefficient fields and $\mathbf{I}_d\in\mathbb{R}^{d\times d}$ is the identity matrix. We assume that $\mathbf{K}(\mathbf{x})$ and $\mathbf{D}(\mathbf{x})$ are symmetric and uniformly positive definite. That is, there exist positive constants $k_{\min}$, $k_{\max}$, $d_{\min}$, and $d_{\max}$ such that
$$
k_{\min}\boldsymbol{\xi}^{T}\boldsymbol{\xi}
\leq
\boldsymbol{\xi}^{T}\mathbf{K}(\mathbf{x})\boldsymbol{\xi}
\leq
k_{\max}\boldsymbol{\xi}^{T}\boldsymbol{\xi},
\qquad
d_{\min}\boldsymbol{\xi}^{T}\boldsymbol{\xi}
\leq
\boldsymbol{\xi}^{T}\mathbf{D}(\mathbf{x})\boldsymbol{\xi}
\leq
d_{\max}\boldsymbol{\xi}^{T}\boldsymbol{\xi},
$$
for all $\boldsymbol{\xi}\in\mathbb{R}^{d}$ and almost every $\mathbf{x}\in\Omega$. In the isotropic case, these conditions reduce to
$$
0<k_{\min}\leq \kappa(\mathbf{x})\leq k_{\max}<\infty,
\qquad
0<d_{\min}\leq d_T(\mathbf{x})\leq d_{\max}<\infty .
$$
Although uniformly bounded, the scalar fields $\kappa$ and $d_T$ may be highly heterogeneous and discontinuous. Representative coefficient fields, displayed on a logarithmic scale, are shown in Figure~\ref{fig:random_KD_fields}.

To complete the model, the system is supplemented with suitable initial and boundary conditions. The initial conditions are prescribed as
$$
\mathbf{u}(\cdot,0)=\mathbf{u}_{0},
\qquad
p(\cdot,0)=p_{0},
\qquad
\theta(\cdot,0)=\theta_{\mathrm{init}}
\qquad
\text{in } \Omega.
$$
The boundary conditions depend on the numerical setting and will be specified in the corresponding numerical experiments.

\subsection{Weak formulation}
\label{subsec:thm_weak_formulation}

We now present the weak formulation of the THM system. Let
\[
\mathbf{V}:=[H^1_0(\Omega)]^d,
\qquad
W:=H^1_0(\Omega),
\qquad
\Theta:=H^1_0(\Omega),
\]
be the function spaces for the displacement, pressure, and temperature variables, respectively. The homogeneous spaces are used here for notational simplicity; the boundary conditions used in the computations are specified in the numerical experiments.

We denote by $(\cdot,\cdot)$ the standard $L^2(\Omega)$ inner product. For scalar functions $a$ and $b$,
\[
(a,b)=\int_{\Omega}ab\,d\mathbf{x}.
\]
For vector- and tensor-valued quantities, the same notation denotes the corresponding componentwise inner product; for example,
\[
(\mathbf{K}\nabla p,\nabla w)
=
\int_{\Omega}\mathbf{K}\nabla p\cdot \nabla w\,d\mathbf{x},
\qquad
(\boldsymbol{\varepsilon}(\mathbf{u}),\boldsymbol{\varepsilon}(\mathbf{v}))
=
\int_{\Omega}
\boldsymbol{\varepsilon}(\mathbf{u}):\boldsymbol{\varepsilon}(\mathbf{v})\,d\mathbf{x}.
\]

Multiplying the momentum, pressure, and temperature equations in \eqref{eq:thm_strong} by test functions
$\mathbf{v}\in\mathbf{V}$, $w\in W$, and $s\in\Theta$, respectively, and integrating over $\Omega$, we obtain the following weak formulation: find
$(\mathbf{u},p,\theta)\in \mathbf{V}\times W\times \Theta$ such that, for all
$(\mathbf{v},w,s)\in \mathbf{V}\times W\times \Theta$,
\begin{subequations}
\label{eq:thm_weak}
\begin{align}
2\mu
\left(
\boldsymbol{\varepsilon}(\mathbf{u}),
\boldsymbol{\varepsilon}(\mathbf{v})
\right)
+
\lambda
\left(
\nabla\cdot\mathbf{u},
\nabla\cdot\mathbf{v}
\right)
-
\alpha
\left(
p,
\nabla\cdot\mathbf{v}
\right)
-
3\alpha_T K_{\mathrm{dr}}
\left(
\theta,
\nabla\cdot\mathbf{v}
\right)
&=
(\mathbf{f},\mathbf{v}),
\label{eq:weak_momentum}
\\
c_0
\left(
p_t,w
\right)
+
\alpha
\left(
\nabla\cdot\mathbf{u}_t,w
\right)
-
3\alpha_m
\left(
\theta_t,w
\right)
+
\left(
\mathbf{K}\nabla p,\nabla w
\right)
&=
(g,w),
\label{eq:weak_pressure}
\\
C_d
\left(
\theta_t,s
\right)
+
3\alpha_T K_{\mathrm{dr}}\theta_0
\left(
\nabla\cdot\mathbf{u}_t,s
\right)
-
3\alpha_m\theta_0
\left(
p_t,s
\right)
+
\left(
\mathbf{D}\nabla\theta,\nabla s
\right)
&=
(\eta,s).
\label{eq:weak_temperature}
\end{align}
\end{subequations}

The first equation represents mechanical equilibrium, where the elastic response is coupled to the fluid pressure and temperature fields. The second equation represents fluid mass conservation, including storage, deformation coupling, thermal dilation, and Darcy diffusion through the permeability field $\mathbf{K}$. The third equation represents energy balance, including heat storage, coupling with deformation and pressure, and heat diffusion through the thermal conductivity field $\mathbf{D}$.

\subsection{Monolithic high-fidelity discretization}
\label{subsec:thm_monolithic_discretization}

We now describe the fully discrete monolithic system used to generate the high-fidelity THM snapshots. Let $\Delta t>0$ be the time step size and let
$$
t^n=n\Delta t,
\qquad
n=0,1,\ldots,N_t,
$$
where $N_t$ denotes the total number of time steps and $t^{N_t}=T$. 
For the spatial discretization, let
$$
\mathbf{V}_h\subset \mathbf{V},
\qquad
W_h\subset W,
\qquad
\Theta_h\subset \Theta
$$
be finite element spaces for the displacement, pressure, and temperature variables, respectively. We use a backward Euler discretization in time. Given
$(\mathbf{u}_h^n,p_h^n,\theta_h^n)\in \mathbf{V}_h\times W_h\times \Theta_h$, the monolithic high-fidelity problem at time $t^{n+1}$ is to find
$(\mathbf{u}_h^{n+1},p_h^{n+1},\theta_h^{n+1})\in \mathbf{V}_h\times W_h\times \Theta_h$ such that, for all
$(\mathbf{v}_h,w_h,s_h)\in \mathbf{V}_h\times W_h\times \Theta_h$,
\begin{subequations}
\label{eq:thm_monolithic_discrete}
\begin{align}
2\mu
\left(
\boldsymbol{\varepsilon}(\mathbf{u}_h^{n+1}),
\boldsymbol{\varepsilon}(\mathbf{v}_h)
\right)
+
\lambda
\left(
\nabla\cdot\mathbf{u}_h^{n+1},
\nabla\cdot\mathbf{v}_h
\right)
-
\alpha
\left(
p_h^{n+1},
\nabla\cdot\mathbf{v}_h
\right)
-
3\alpha_T K_{\mathrm{dr}}
\left(
\theta_h^{n+1},
\nabla\cdot\mathbf{v}_h
\right)
&=
(\mathbf{f}^{n+1},\mathbf{v}_h),
\label{eq:discrete_momentum}
\\
c_0
\left(
\frac{p_h^{n+1}-p_h^n}{\Delta t},
w_h
\right)
+
\alpha
\left(
\frac{\nabla\cdot\mathbf{u}_h^{n+1}-\nabla\cdot\mathbf{u}_h^n}{\Delta t},
w_h
\right)
-
3\alpha_m
\left(
\frac{\theta_h^{n+1}-\theta_h^n}{\Delta t},
w_h
\right)
+
\left(
\mathbf{K}\nabla p_h^{n+1},
\nabla w_h
\right)
&=
(g^{n+1},w_h),
\label{eq:discrete_pressure}
\\
C_d
\left(
\frac{\theta_h^{n+1}-\theta_h^n}{\Delta t},
s_h
\right)
+
3\alpha_TK_{\mathrm{dr}}\theta_0
\left(
\frac{\nabla\cdot\mathbf{u}_h^{n+1}-\nabla\cdot\mathbf{u}_h^n}{\Delta t},
s_h
\right)
-
3\alpha_m\theta_0
\left(
\frac{p_h^{n+1}-p_h^n}{\Delta t},
s_h
\right)
+
\left(
\mathbf{D}\nabla\theta_h^{n+1},
\nabla s_h
\right)
&=
(\eta^{n+1},s_h).
\label{eq:discrete_temperature}
\end{align}
\end{subequations}

The system \eqref{eq:thm_monolithic_discrete} is monolithic because the displacement, pressure, and temperature unknowns at time $t^{n+1}$ are solved simultaneously. To write the fully discrete system in algebraic form, let
$\{\boldsymbol{\varphi}_i^u\}_{i=1}^{N_u}$,
$\{\varphi_i^p\}_{i=1}^{N_p}$, and
$\{\varphi_i^\theta\}_{i=1}^{N_\theta}$
be finite element bases for $\mathbf{V}_h$, $W_h$, and $\Theta_h$, respectively. Here,
$N_u$, $N_p$, and $N_\theta$ denote the numbers of degrees of freedom for the displacement, pressure, and temperature variables. At time $t^n$, we write
$$
\mathbf{u}_h^n
=
\sum_{j=1}^{N_u} U_j^n \boldsymbol{\varphi}_j^u,
\qquad
p_h^n
=
\sum_{j=1}^{N_p} P_j^n \varphi_j^p,
\qquad
\theta_h^n
=
\sum_{j=1}^{N_\theta} \Theta_j^n \varphi_j^\theta .
$$
The scalars $U_j^n$, $P_j^n$, and $\Theta_j^n$ are the finite element coefficients of the displacement, pressure, and temperature approximations at time $t^n$. We collect them into the vectors
$$
\mathbf{U}^n
=
\begin{bmatrix}
U_1^n & U_2^n & \cdots & U_{N_u}^n
\end{bmatrix}^{T}
\in \mathbb{R}^{N_u},
$$
$$
\mathbf{P}^n
=
\begin{bmatrix}
P_1^n & P_2^n & \cdots & P_{N_p}^n
\end{bmatrix}^{T}
\in \mathbb{R}^{N_p},
\qquad
\boldsymbol{\Theta}^n
=
\begin{bmatrix}
\Theta_1^n & \Theta_2^n & \cdots & \Theta_{N_\theta}^n
\end{bmatrix}^{T}
\in \mathbb{R}^{N_\theta}.
$$
The monolithic state vector at time $t^n$ is then defined by
$$
\mathbf{z}_h^n
=
\begin{bmatrix}
\mathbf{U}^n\\
\mathbf{P}^n\\
\boldsymbol{\Theta}^n
\end{bmatrix}
\in \mathbb{R}^{N_h},
\qquad
N_h=N_u+N_p+N_\theta .
$$

With these definitions, the fully discrete THM system can be written as the block linear system
\begin{equation}
\label{eq:monolithic_block_system}
\begin{bmatrix}
\mathbf{A}_{uu} & \mathbf{A}_{up} & \mathbf{A}_{u\theta}
\\
\mathbf{M}_{pu} & \mathbf{M}_{pp}+\mathbf{A}_{pp} & \mathbf{M}_{p\theta}
\\
\mathbf{M}_{\theta u} & \mathbf{M}_{\theta p} & \mathbf{M}_{\theta\theta}+\mathbf{A}_{\theta\theta}
\end{bmatrix}
\begin{bmatrix}
\mathbf{U}^{n+1}\\
\mathbf{P}^{n+1}\\
\boldsymbol{\Theta}^{n+1}
\end{bmatrix}
=
\begin{bmatrix}
\mathbf{F}^{n+1}
\\
\mathbf{G}^{n+1}
+
\mathbf{M}_{pu}\mathbf{U}^{n}
+
\mathbf{M}_{pp}\mathbf{P}^{n}
+
\mathbf{M}_{p\theta}\boldsymbol{\Theta}^{n}
\\
\mathbf{Q}^{n+1}
+
\mathbf{M}_{\theta u}\mathbf{U}^{n}
+
\mathbf{M}_{\theta p}\mathbf{P}^{n}
+
\mathbf{M}_{\theta\theta}\boldsymbol{\Theta}^{n}
\end{bmatrix}.
\end{equation}
Here, $\mathbf{A}_{uu}\in\mathbb{R}^{N_u\times N_u}$ is the elasticity matrix associated with the displacement unknown. The matrices
$\mathbf{A}_{up}\in\mathbb{R}^{N_u\times N_p}$ and
$\mathbf{A}_{u\theta}\in\mathbb{R}^{N_u\times N_\theta}$
represent the coupling of the mechanics equation with the pressure and temperature variables, respectively. The matrices
$\mathbf{M}_{pp},\mathbf{A}_{pp}\in\mathbb{R}^{N_p\times N_p}$
are the pressure storage and pressure diffusion matrices, while
$\mathbf{M}_{pu}\in\mathbb{R}^{N_p\times N_u}$ and
$\mathbf{M}_{p\theta}\in\mathbb{R}^{N_p\times N_\theta}$
are the coupling blocks in the pressure equation associated with displacement and temperature. Similarly,
$\mathbf{M}_{\theta\theta},\mathbf{A}_{\theta\theta}\in\mathbb{R}^{N_\theta\times N_\theta}$
are the temperature storage and heat diffusion matrices, while
$\mathbf{M}_{\theta u}\in\mathbb{R}^{N_\theta\times N_u}$ and
$\mathbf{M}_{\theta p}\in\mathbb{R}^{N_\theta\times N_p}$
are the coupling blocks in the temperature equation associated with displacement and pressure. Finally,
$\mathbf{F}^{n+1}\in\mathbb{R}^{N_u}$,
$\mathbf{G}^{n+1}\in\mathbb{R}^{N_p}$, and
$\mathbf{Q}^{n+1}\in\mathbb{R}^{N_\theta}$
denote the discrete body force, fluid source/sink, and heat source/sink vectors, respectively.

Equivalently, the block system \eqref{eq:monolithic_block_system} can be written in compact form as
$$
\mathcal{A}_h \mathbf{z}_h^{n+1}
=
\mathbf{b}_h^{n+1}(\mathbf{z}_h^n),
$$
where $\mathcal{A}_h\in\mathbb{R}^{N_h\times N_h}$ is the fully coupled monolithic THM system matrix, $\mathbf{z}_h^{n+1}\in\mathbb{R}^{N_h}$ is the monolithic state vector at time $t^{n+1}$, and $\mathbf{b}_h^{n+1}(\mathbf{z}_h^n)\in\mathbb{R}^{N_h}$ contains the source terms together with the contributions from the previous time-step solution $\mathbf{z}_h^n\in\mathbb{R}^{N_h}$.

Solving this system for $n=0,1,\ldots,N_t-1$ produces the high-fidelity monolithic solution sequence
$$
\mathbf{z}_h^1,\mathbf{z}_h^2,\ldots,\mathbf{z}_h^{N_t},
\qquad
\mathbf{z}_h^n\in\mathbb{R}^{N_h}.
$$

These solution vectors are then assembled into the snapshot matrix
$$
\mathbf{Z}
=
\begin{bmatrix}
\mathbf{z}_h^1 & \mathbf{z}_h^2 & \cdots & \mathbf{z}_h^{N_t}
\end{bmatrix}
\in \mathbb{R}^{N_h\times N_t}.
$$
The matrix $\mathbf{Z}$ is the high-fidelity data matrix used as the input for SparseModesNet in the next section.

 \section{SparseModesNet Approximation for THM}\label{sec:SparseModesNet} 

SparseModesNet, originally introduced in \cite{koike2026SparseModesNet}, is a dimensionality-reduction method that combines a linear POD encoder with a sparsifying nonlinear NN decoder. Unlike standard POD, SMN does not necessarily use the first $r$ modes. Instead, the method first constructs a larger set of $s$ candidate POD modes and then select $r<s$ modes that are most useful for reconstructing the solution data. Importantly, the NN and selected POD modes are optimized simultaneously.

Starting from the snapshot matrix $\mathbf{Z}
=
\begin{bmatrix}
\mathbf{z}_h^1 & \mathbf{z}_h^2 & \cdots & \mathbf{z}_h^{N_t}
\end{bmatrix}
\in \mathbb{R}^{N_h\times N_t},$ defined in section \ref{sec:thm_model}. We compute its rank-s truncated singular value decomposition  
\begin{equation}
    \mathbf{Z}
\approx
\mathbf{U}_{s}\boldsymbol{\Sigma}_{s}\mathbf{V}_{s}^{\top},
\end{equation}
where, $\mathbf{U}_{s}\in\mathbb{R}^{N_h\times s}$, $\boldsymbol{\Sigma}_{s}\in\mathbb{R}^{s\times s},$ 
and $\mathbf{V}_{s}\in\mathbb{R}^{N_t\times s}.$ The columns of $
\mathbf{U}_{s}
=
\begin{bmatrix}
\boldsymbol{\phi}_{0} &
\boldsymbol{\phi}_{1} &
\cdots &
\boldsymbol{\phi}_{s-1}
\end{bmatrix}
$, are the first $s$ POD modes and form the candidate set used by SMN. For each snapshot $\mathbf{z}_h^n$, the corresponding coordinates in this candidate basis are 
\begin{equation}
    \widehat{\mathbf{a}}^{\,n}
=
\mathbf{U}_{s}^{\top}\mathbf{z}_h^n
\in\mathbb{R}^{s}.
\end{equation}

For comparison, standard POD retains the first $r$ modes, $\mathbf{U}_r^{\mathrm{POD}}
=
\begin{bmatrix}
\boldsymbol{\phi}_0 & \cdots & \boldsymbol{\phi}_{r-1}
\end{bmatrix} \in \mathbb{R}^{N_h \times r},$ and we denote the POD reconstruction of the snapshot \(\mathbf{z}_h^n\) by
\(\widetilde{\mathbf{z}}_{\mathrm{POD}}^{\,n}\), which is given by

 \begin{equation}
\widetilde{\mathbf{z}}_{\mathrm{POD}}^{\,n}
=
\mathbf{U}_r^{\mathrm{POD}}
\bigl(\mathbf{U}_r^{\mathrm{POD}}\bigr)^{\top}
\mathbf{z}_h^n.
\end{equation}
SMN instead selects \(r\) modes from the larger set of \(s\)
candidate modes.

\subsection{Sparse mode selection}\label{subsec:mode_selection}

SMN select $r$ modes from the $s$ candidate modes by introducing a trainable mode-selection vector $\boldsymbol{\omega}
=
\begin{bmatrix}
\omega_0 & \cdots & \omega_{s-1}
\end{bmatrix}^{\top}
\in\mathbb{R}^{s}$. This vector acts componentwise on the candidate coordinates:
$$
\boldsymbol{\omega}\odot\widehat{\mathbf{a}}^{\,n},
$$
where $\odot$ denotes component wise multiplication. If $\omega_j=0$, then the coordinate associated with the $j$-th candidate POD mode is inactive.

During mode selection, SMN decoder
\begin{equation}
\label{eq:smn-selection-decoder}
\mathcal{D}_{\mathrm{sel}}
\left(
\widehat{\mathbf{a}}^{\,n}
\right)
=
\underbrace{
\mathbf{U}_s
\left(
\boldsymbol{\omega}
\odot
\widehat{\mathbf{a}}^{\,n}
\right)
}_{\text{linear POD contribution}}
+ 
\underbrace{
\mathbf{W}_{\mathrm{sel}}
h_{\boldsymbol{\vartheta}}
\left(
\boldsymbol{\omega}
\odot
\widehat{\mathbf{a}}^{\,n}
\right)
}_{\text{nonlinear correction}}.
\end{equation}
Here, $h_{\boldsymbol{\vartheta}}:
\mathbb{R}^{s}\rightarrow\mathbb{R}^{q}$ is a neural network with parameter $\boldsymbol{\vartheta}$, and $\mathbf{W}_{\mathrm{sel}}
\in\mathbb{R}^{N_h\times q}$ maps the neural-network output to the original finite element solution space
\(\mathbb{R}^{N_h}\). The same weighted coordinates are used in both parts of the decoder. Therefore, if $\omega_j=0$, the $j$-th candidate coordinate cannot contribute to either the linear reconstruction or the nonlinear correction. For a given $\omega$, define the active mode set by 
\begin{equation}
\label{eq:smn-active-set}
\mathcal{A}(\boldsymbol{\omega})
=
\left\{
j\in\{0,\ldots,{s-1}\}:\omega_j\neq0
\right\}.
\end{equation}
We denote by $\mathbf{U}_{\mathcal{A}(\boldsymbol{\omega})}$ the matrix containing the column of $\mathbf{U}_s$, whose indices belong to $\mathcal{A}(\boldsymbol{\omega})$.

Let $\boldsymbol{\Theta}^{(0)}$ denote the first-layer weight matrix of $h_{\boldsymbol{\vartheta}}$, and let $\boldsymbol{\Theta}_{:,j}^{(1)}$ denote the column containing all first-layer weights connected to the $j$-th POD coordinate. Following the LassoNet \cite{lemhadri2021lassonet} formulation, the mode-selection vector and the nonlinear decoder are trained jointly by solving 
\begin{subequations}
\label{eq:smn-mode-selection}
\begin{align}
\min_{
\boldsymbol{\omega},
\boldsymbol{\vartheta},
\mathbf{W}_{\mathrm{sel}}
}
\quad&
\frac{1}{N_t}
\sum_{n=1}^{N_t}
\left\|
\mathbf{z}_h^n
-
\mathcal{D}_{\mathrm{sel}}
\left(
\widehat{\mathbf{a}}^{\,n}
\right)
\right\|_2^2
+
\lambda\|\boldsymbol{\omega}\|_1,
\label{eq:smn-selection-objective}
\\
\text{subject to}\quad&
\left\|
\boldsymbol{\Theta}_{:,j}^{(0)}
\right\|_{\infty}
\leq
M|\omega_j|,
\qquad j=0,\ldots,s-1,
\label{eq:smn-hierarchical-constraint}
\\
&
\mathbf{U}_{\mathcal{A}(\boldsymbol{\omega})}^{\top}
\mathbf{W}_{\mathrm{sel}}
=
\mathbf{0}.
\label{eq:smn-orthogonality-constraint}
\end{align}
\end{subequations}
Here, $\lambda\geq0$ controls the sparsity of $\boldsymbol{\omega}$, while $M>0$ is the parameters.  The $\ell_1$-penalty in~\eqref{eq:smn-selection-objective}  makes some entries of $\boldsymbol{\omega}$ to become zero. If $\omega_j=0$, the constraint in~\eqref{eq:smn-hierarchical-constraint} also forces all first-layer weight is connected to the $j$-th POD coordinate to be zero, i.e,  $\omega_j=0
~\text{then}~
\boldsymbol{\Theta}_{:,j}^{(0)}=\mathbf{0}.$ Therefore, the $j$-th POD coordinate is removed from both the linear and nonlinear parts of the decoder. The orthogonality constraint in \eqref{eq:smn-orthogonality-constraint}  ensures that the nonlinear correction represents information that is not already contained in the active POD modes. 

Problem~\eqref{eq:smn-selection-objective} is solved for a sequence of increasing values of $\lambda$. Starting from $\lambda_0$, we update 
\begin{equation}
\label{eq:smn-continuation}
\lambda
\leftarrow
(1+\varepsilon)\lambda,
\qquad
\varepsilon>0.
\end{equation}
For each new value of $\lambda$, training starts from the solution obtained for the previous value.  
As $\lambda$ increases, more entries of $\boldsymbol{\omega}$ become zero and the corresponding POD modes are removed. This process continues until $r$ modes remain active. 

Let $\boldsymbol{\omega}^*$ denote the selection vector obtained at the end of this procedure. The indices of $r$ active modes are 
\begin{equation}
\label{eq:smn-selected-indices}
\mathcal{I}_{SMN}
=
\left\{
j\in\{0,\ldots,s-1\}:
\omega_j^*\neq0
\right\},
\qquad
|\mathcal{I}_{SMN}|=r.
\end{equation}
If $\mathcal{I}_{SMN}=\{j_0,\ldots,j_{r-1}\},$ then the selected POD basis is
\begin{equation}
\label{eq:smn-selected-basis}
\mathbf{U}_r
=
\begin{bmatrix}
\boldsymbol{\phi}_{j_0} &
\cdots &
\boldsymbol{\phi}_{j_{r-1}}
\end{bmatrix}
\in\mathbb{R}^{N_h\times r}.
\end{equation}
The columns of $\mathbf{U}_r$ are therefore the $r$ POD modes selected from the $s$ candidate modes. 

\subsection{Nonlinear decoder retraining} 

The parameters learned during the selection phase are used only to determine the selected mode set $\mathcal{I}_r$. After fixing $\mathbf{U}_r$, SparseModesNet retrains the nonlinear decoder from scratch while dropping the constraint in equation \eqref{eq:smn-hierarchical-constraint} and $\ell_1$- penalty to avoid the slight downward bias introduced by the regularization~\cite{lemhadri2021lassonet}.

\begin{figure}[!h]
\centering
\includegraphics[width=\textwidth]
{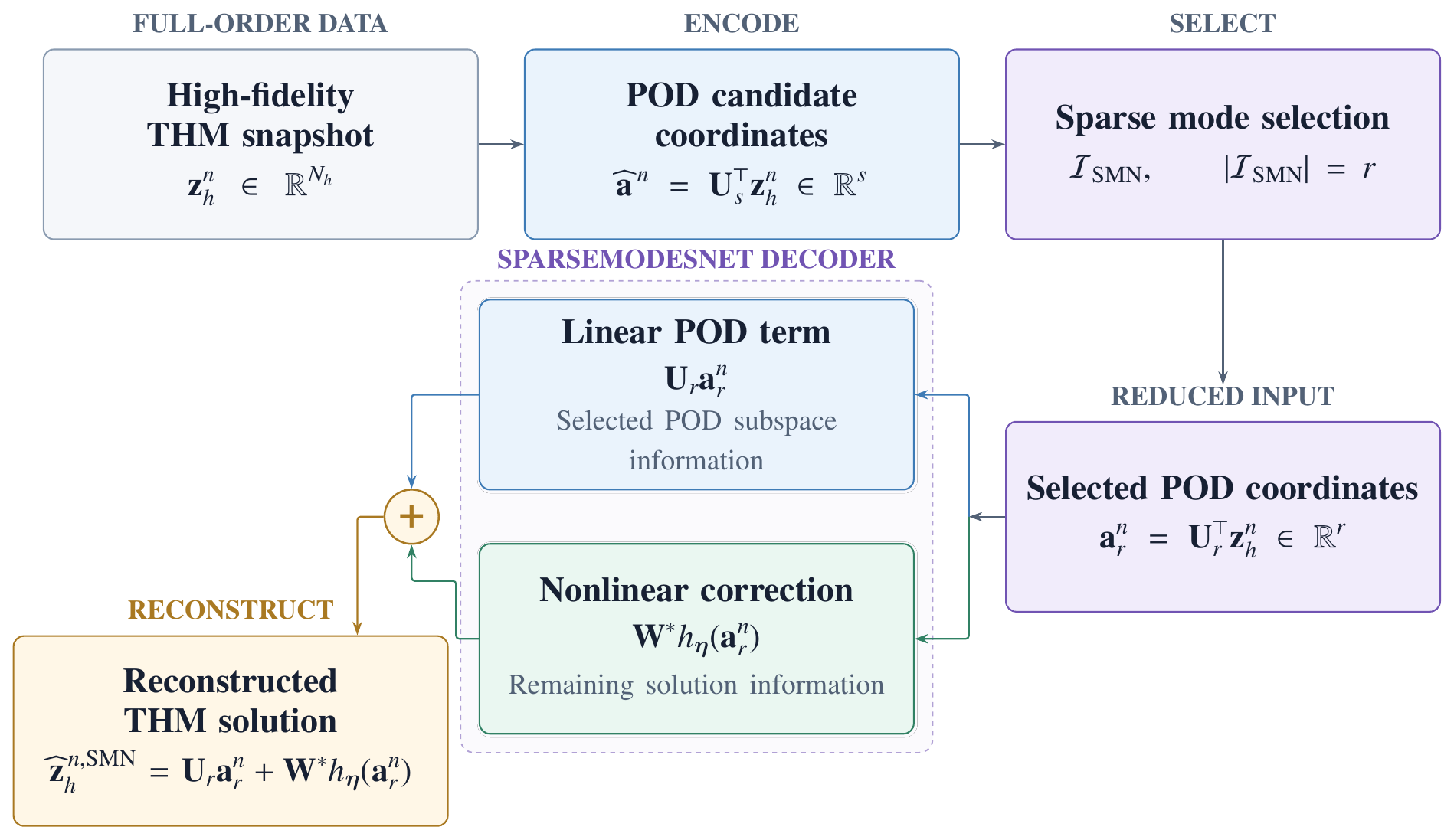}
\caption{Workflow of the SparseModesNet approximation. The selected POD
coordinates are passed to the linear POD term and nonlinear correction, whose
sum gives the reconstructed THM state.}
\label{fig:smn-workflow}
\end{figure}

For each snapshot, the selected reduced coordinates are \begin{equation}
\label{eq:smn-selected-coordinates}
\mathbf{a}_r^n
=
\mathbf{U}_r^{\top}\mathbf{z}_h^n
\in\mathbb{R}^{r}.
\end{equation}
The final decoder has the form 
\begin{equation}
\label{eq:smn-final-decoder}
\mathcal{D}_{\mathrm{SMN}}(\mathbf{a}_r^n)
=
\mathbf{U}_r\mathbf{a}_r^n
+
\mathbf{W}
h_{\boldsymbol{\eta}}(\mathbf{a}_r^n),
\end{equation}
where, $h_{\boldsymbol{\eta}}:
\mathbb{R}^{r}\rightarrow\mathbb{R}^{q}$ is retrained neural network and $\mathbf{W}\in\mathbb{R}^{N_h\times q}.$ The decoder parameters are determined from 
\begin{equation}
\label{eq:smn-decoder-training}
\begin{aligned}
\min_{\boldsymbol{\eta},\,\mathbf{W}}
\quad&
\frac{1}{N_t}
\sum_{n=1}^{N_t}
\left\|
\mathbf{z}_h^n
-
\mathbf{U}_r\mathbf{a}_r^n
-
\mathbf{W}
h_{\boldsymbol{\eta}}(\mathbf{a}_r^n)
\right\|_2^2,
\\
\text{subject to}\quad&
\mathbf{U}_r^{\top}\mathbf{W}
=
\mathbf{0}.
\end{aligned}
\end{equation}
After training $h_{\boldsymbol{\eta}}$, we fix the trained NN and then the final output matrix is computed from the regularized least-squares problem 
\begin{equation}
\label{eq:smn-output-matrix}
\mathbf{W}^*
=
\arg\min_{\mathbf{W}}
\left\|
\mathbf{Z}
-
\mathbf{U}_r\mathbf{U}_r^{\top}\mathbf{Z}
-
\mathbf{W}
h_{\boldsymbol{\eta}}
\left(
\mathbf{U}_r^{\top}\mathbf{Z}
\right)
\right\|_F^2
+
\gamma\|\mathbf{W}\|_F^2,
\end{equation}
where $\gamma\geq0$ is a regularization parameter. The final SparseModesNet reconstruction of the $n$-th snapshot is therefore 
\begin{equation}
\label{eq:smn-final-reconstruction}
\widehat{\mathbf{z}}_h^{\,n,\mathrm{SMN}}
=
\mathbf{U}_r\mathbf{a}_r^n
+
\mathbf{W}^*
h_{\boldsymbol{\eta}}(\mathbf{a}_r^n).
\end{equation}
The linear term represents the part of the solution contained in the selected POD subspace, while the nonlinear term learns the remaining information needed to reconstruct the complete high-fidelity THM solution. Figure~\ref{fig:smn-workflow} summarizes the SMN workflow.

\section{Numerical Experiments}
\label{sec:numerical}

In this section, we evaluate SMN using solution snapshots generated
by the monolithic high-fidelity THM solver described in section~\ref{sec:thm_model}. The high-fidelity systems are
assembled and solved using FEniCS~\cite{logg2012automated}, while SparseModesNet is implemented in
PyTorch~\cite{paszke2019pytorch}.

Unless stated otherwise, all numerical experiments use homogeneous initial conditions, $\mathbf{u}(\mathbf{x},0)=\mathbf{0}, p(\mathbf{x},0)=0,$ and $\theta(\mathbf{x},0)=0.$ The displacement is fixed on the boundary, $\mathbf{u}=\mathbf{0}$ on \(\partial\Omega\times \mathcal{J}\) , while the
homogeneous no-flux conditions
\(\mathbf{K}\nabla p\cdot\mathbf{n}=0\) and
\(\mathbf{D}\nabla\theta\cdot\mathbf{n}=0\) are imposed for pressure and
temperature on \(\partial\Omega\times \mathcal{J}\). Here, \(\mathbf{n}\) denotes the
outward unit normal vector. The body force in the mechanical equation is set
to zero, \(\mathbf{f}(\mathbf{x},t)=\mathbf{0}\). The fixed THM parameters used in this numerical experiments are 
$
\lambda=\mu=10^2,
K_{\mathrm{dr}}=2\times10^2,
\alpha=c_0=C_d=\theta_0=1,
\alpha_T=10^{-3},$ and 
$\alpha_m=10^{-5}.$
The computational domain, coefficient fields, source term, spatial mesh and time discretization are specified separately for each experiment. 

For each experiment, the computed finite element solutions are collected into a snapshot matrix. Using the notation from Section~\ref{sec:thm_model}, the monolithic state at \(t^n\) and
the corresponding snapshot matrix are written as
\[
\mathbf{z}_h^n
=
\begin{bmatrix}
\mathbf{U}^n\\
\mathbf{P}^n\\
\boldsymbol{\Theta}^n
\end{bmatrix}
\in\mathbb{R}^{N_h},
\qquad
\mathbf{Z}
=
\begin{bmatrix}
\mathbf{z}_h^1&\cdots&\mathbf{z}_h^{N_t}
\end{bmatrix}
\in\mathbb{R}^{N_h\times N_t}.
\]
Here, \(\mathbf{U}^n\), \(\mathbf{P}^n\), and
\(\boldsymbol{\Theta}^n\) are the finite element coefficient vectors
associated with the finite element solutions
\(\mathbf{u}_h(\mathbf{x},t^n)\), \(p_h(\mathbf{x},t^n)\), and
\(\theta_h(\mathbf{x},t^n)\), respectively.

To compare POD and SMN, let \(\mathcal{M}\in\{\mathrm{POD},\mathrm{SMN}\}\) denote the reconstruction
method, and let \(\widetilde{\mathbf{Z}}^{\mathcal{M}}\) denote the
corresponding reconstructed full-state snapshot matrix. The relative reconstruction error over the complete snapshot matrix is defined by
\begin{equation}
    E^{\mathcal{M}}
=
\frac{
\|\mathbf{Z}-\widetilde{\mathbf{Z}}^{\mathcal{M}}\|_{\mathrm{F}}
}{
\|\mathbf{Z}\|_{\mathrm{F}}
},
\end{equation}
where $\|\cdot\|_{\mathrm{F}}$ denotes the Frobenius norm.

We also compute reconstruction error for each solution component. For a component $q$, let $\mathbf{Z}_q=[\mathbf{q}_h^1,\ldots,\mathbf{q}_h^{N_t}]$ and $\widetilde{\mathbf{Z}}_q^{\mathcal{M}}$ denotes its high-fidelity and reconstructed snapshot matrices, respectively. The corresponding component-wise error is 
\begin{equation}
    E_q^{\mathcal{M}}
=
\frac{
\|\mathbf{Z}_q-\widetilde{\mathbf{Z}}_q^{\mathcal{M}}\|_{\mathrm{F}}
}{
\|\mathbf{Z}_q\|_{\mathrm{F}}
}.
\end{equation}

To evaluate the reconstruction at an individual time step, we also use the
pointwise absolute error and the snapshot relative error. For
$q\in\{p,\theta,|\mathbf{u}|\}$, these errors are defined as
\begin{equation}
\label{eq:snapshot-error-measures}
e_q^{\mathcal{M}}(\mathbf{x},t^n)
=
\left|
q_h(\mathbf{x},t^n)
-
\widetilde{q}_h^{\mathcal{M}}(\mathbf{x},t^n)
\right|,
\qquad
\varepsilon_q^{\mathcal{M}}(t^n)
=
\frac{
\left\|
\mathbf{q}_h^n
-
\widetilde{\mathbf{q}}_h^{\mathcal{M},n}
\right\|_2
}{
\left\|
\mathbf{q}_h^n
\right\|_2
},
\end{equation}
where $\mathcal{M}\in\{\mathrm{POD},\mathrm{SMN}\}$.
Here, $\mathbf{q}_h^n$ and
$\widetilde{\mathbf{q}}_h^{\mathcal{M},n}$ denote the finite element
coefficient vectors of the high-fidelity solution and its reconstruction at
time $t^n$, respectively. The notation $\|\cdot\|_2$ denotes the Euclidean
norm. 
The following subsections present the numerical experiments and compare POD and SMN using these error measure. 

\subsection{Experiment 1: THM system with heterogeneous permeability and thermal conductivity}
\label{subsec:exp1_thm_heterogeneous}

\subsubsection{Experiment setup and time-dependent well operations}

The first experiment is a synthetic benchmark constructed to test whether SMN can reconstruct THM solutions with randomly generated, spatially heterogeneous permeability and thermal conductivity fields, together with time-dependent well operations. The design of this benchmark is motivated by geothermal studies of coupled THM
behavior during injection and production, reinjection allocation, time-dependent
injection scheduling, and interactions between injection and production
wells~\cite{akin2014reinjection,juliusson2013optimization,
daniilidis2021interference,lei2023injection}.

\begin{figure}[!h]
\centering
\begin{subfigure}[t]{0.32\textwidth}
\centering
\includegraphics[width=\linewidth]{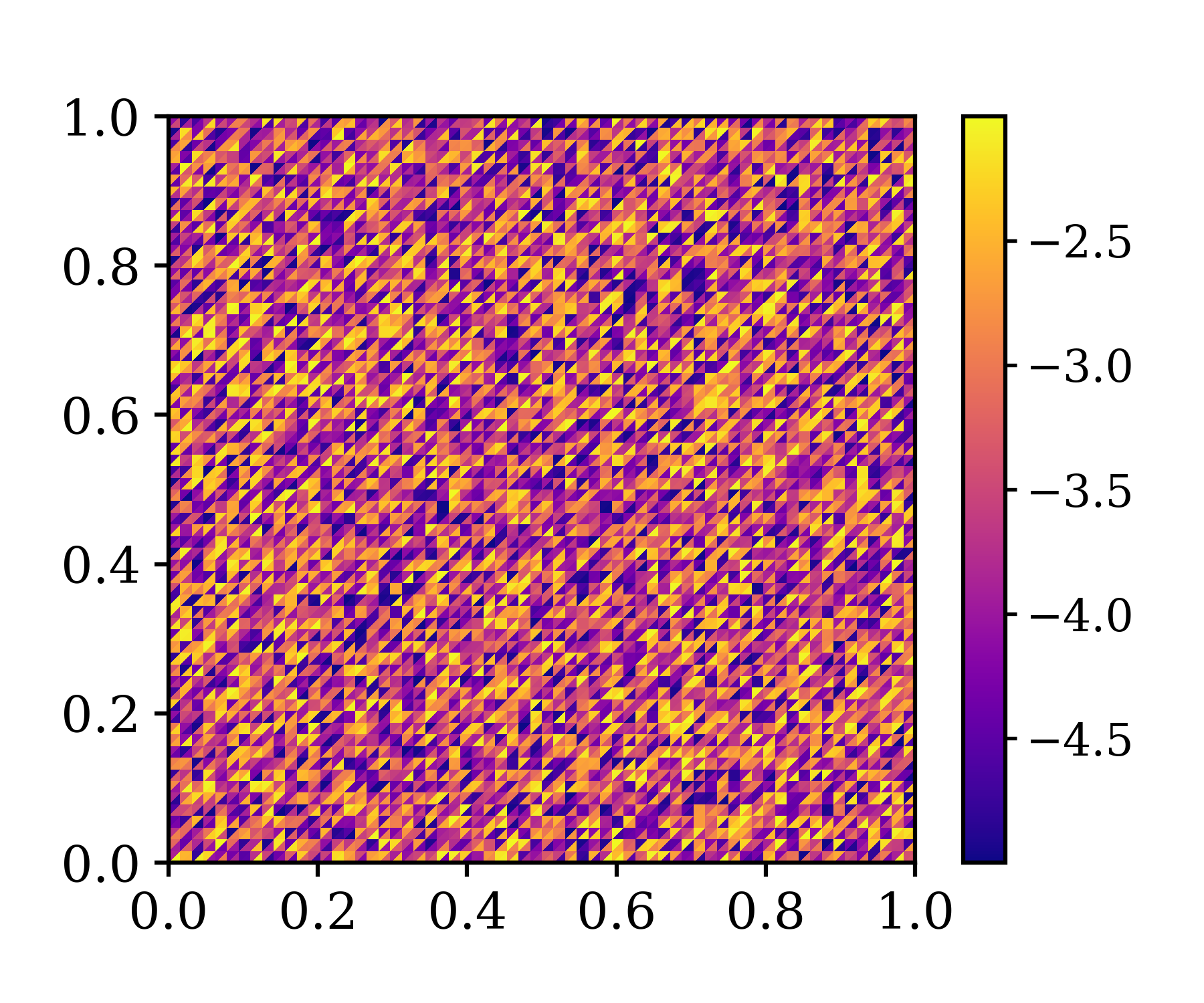}
\caption{Permeability \(\kappa(\mathbf{x})\).}
\label{fig:exp1-perm}
\end{subfigure}
\hfill
\begin{subfigure}[t]{0.32\textwidth}
\centering
\includegraphics[width=\linewidth]{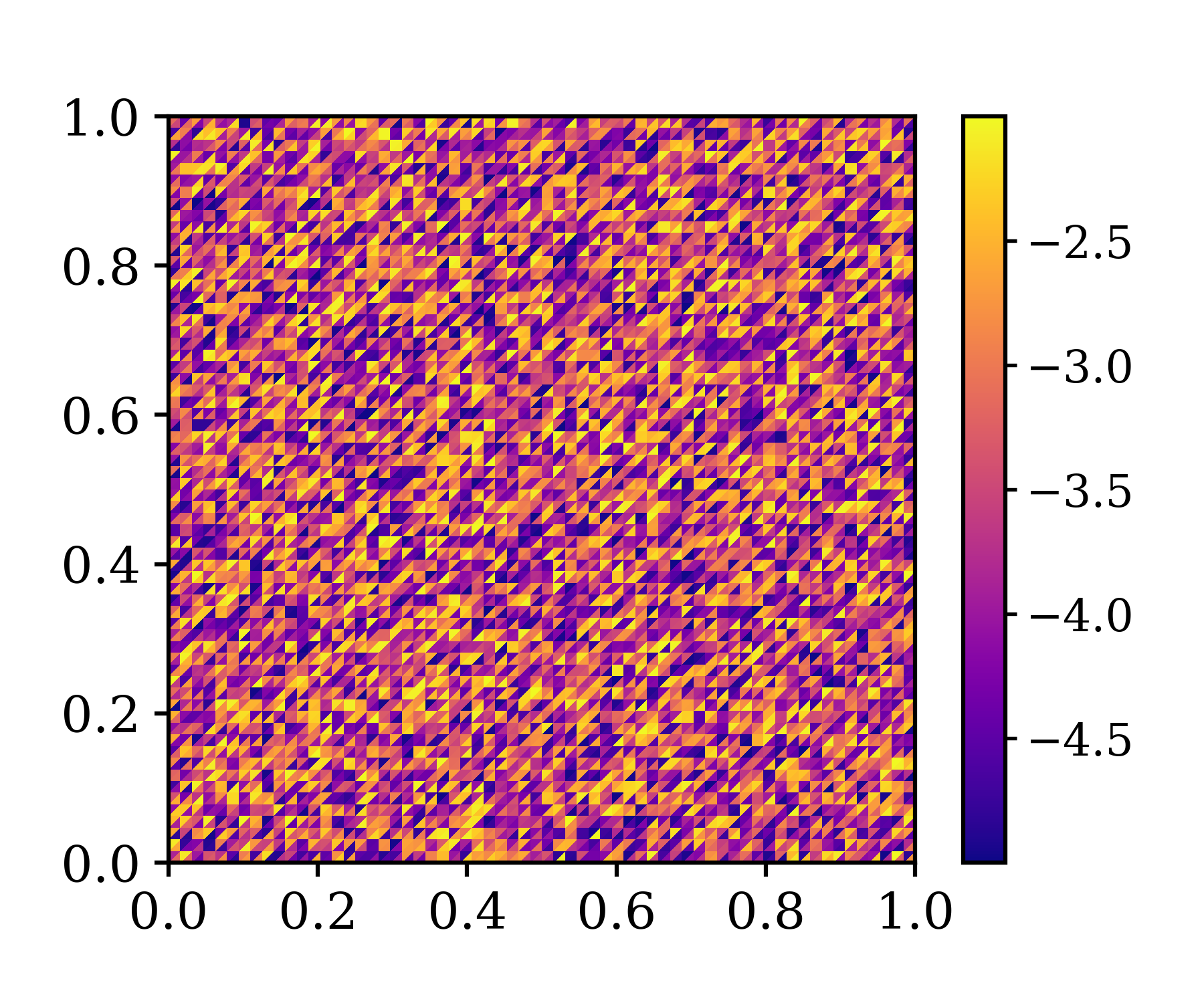}
\caption{Thermal conductivity \(d_T(\mathbf{x})\).}
\label{fig:exp1-thermal}
\end{subfigure}
\hfill
\begin{subfigure}[t]{0.32\textwidth}
\centering
\includegraphics[width=\linewidth]{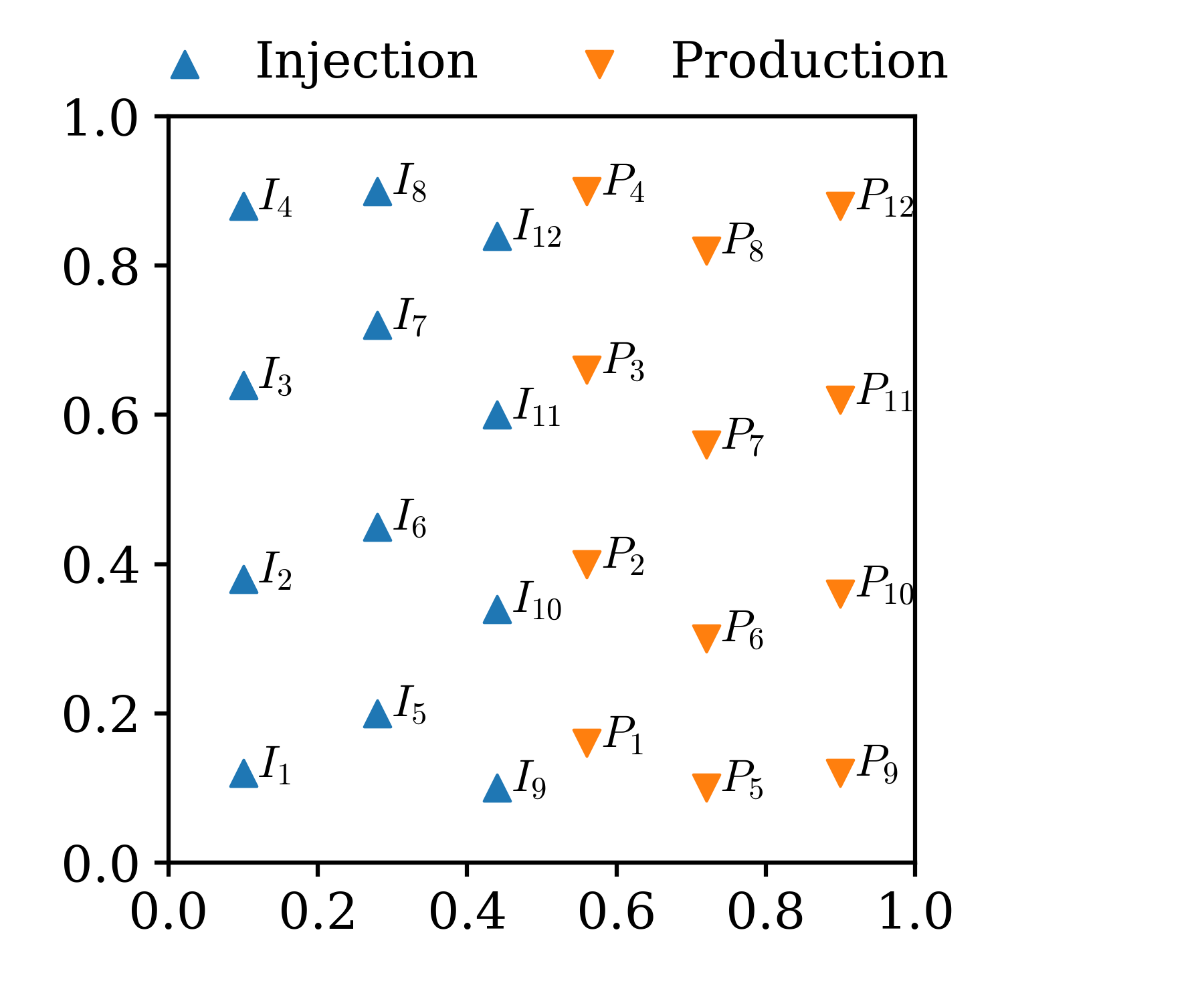}
\caption{Fixed injection-production locations.}
\label{fig:exp1-well-locations}
\end{subfigure}
\caption{Experiment~1. Heterogeneous permeability and thermal conductivity
fields displayed on a logarithmic scale, together with the fixed injection
and production locations.}
\label{fig:thm-random-kd-coefficients}
\end{figure}

We solve the same THM problem described in Section \ref{sec:thm_model} on $\Omega=(0,1)^2$ over $\mathcal{J}=(0,2]$. Permeability and thermal conductivity  coefficients are isotropic, with $\mathbf{K}(\mathbf{x})=\kappa(\mathbf{x})\mathbf{I}_2$ and $\mathbf{D}(\mathbf{x})=d_T(\mathbf{x})\mathbf{I}_2$, where $\mathbf{I}_2$ is a $2\times2$ identity matrix. The values of $\kappa(\mathbf{x})$ and $d_T(\mathbf{x})$ vary randomly from one mesh element to another. Figures~\ref{fig:exp1-perm} and \ref{fig:exp1-thermal} show the
coefficient fields used in this experiment.

In this experiment,  we use $12$ fixed injection locations and $12$ production locations, as shown in
Figure~\ref{fig:exp1-well-locations}. Let $\mathcal{C}^{\mathrm{I}}$ and $\mathcal{C}^{\mathrm{P}}$ denote the set of injection and production coordinates, respectively. We introduce a
time-dependent operating schedule in which different injection-production
pairs are activated at different times. This changes the locations and
magnitudes of the source terms throughout the simulation. As a result, the
pressure, temperature, and displacement fields exhibit spatially heterogeneous
and time-dependent behavior, producing a challenging data set for
reconstruction.

We divide the time interval $\mathcal{J}=(0,2]$ into $20$ operating stages, each of length $0.1$. During stage $m$, let $\mathbf{c}_m^{\mathrm{I}}\in\mathcal{C}^{\mathrm{I}}$ and $\mathbf{c}_m^{\mathrm{P}}\in\mathcal{C}^{\mathrm{P}}$ denote the active injection and production locations, respectively. The target rate for this stage is $Q_m=a_mQ_{\text{ref}}$, where $
Q_{\text{ref}}=1.8\times10^{-2},$ and $a_m$ is the rate factor assigned to stage $m$. Let $t_m=(m-1)\Delta T$ denotes the beginning of stage $m$. For $m>2$, the operation is transferred from the previous pair $(\mathbf{c}_{m-1}^{\mathrm{I}},\mathbf{c}_{m-1}^{\mathrm{P}})$ to new pair $(\mathbf{c}_m^{\mathrm{I}},\mathbf{c}_m^{\mathrm{P}})$ over the transition interval $[t_m,t_m+\delta]$, where $\delta=0.02$. For
\(t\in[t_m,t_m+\delta]\), we set
\(\tau=(t-t_m)/\delta\) and \(\xi(\tau)=3\tau^2-2\tau^3\). The previous
and new pairs operate with rates
\(Q_{\mathrm{old}}(t)=(1-\xi(\tau))Q_{m-1}\) and
\(Q_{\mathrm{new}}(t)=\xi(\tau)Q_m\), respectively. Thus, the common total injection and production rate during the transition is $ Q(t) = Q_{\mathrm{old}}(t)+Q_{\mathrm{new}}(t).$
The previous and new pairs receive \(Q_{\mathrm{old}}(t)\) and \(Q_{\mathrm{new}}(t)\), respectively, while all other locations have zero rate. After the transition, the previous pair is turned off and the new pair operates with rate \(Q_m\) for the remainder of the stage. In the first stage, the rate is increased smoothly from zero to \(Q_1\).

Each location $\mathbf{c}\in\mathcal{C}^{\mathrm{I}}\cup\mathcal{C}^{\mathrm{P}}$ is represented by normalized Gaussian function 
\begin{equation}\label{eq:gaussian_fun}
G(\mathbf{x};\mathbf{c})
=
\frac{
\exp\!\left(
-\beta_{\mathrm{src}}\|\mathbf{x}-\mathbf{c}\|_2^2
\right)
}{
\displaystyle
\int_{\Omega}
\exp\!\left(
-\beta_{\mathrm{src}}\|\mathbf{y}-\mathbf{c}\|_2^2
\right)\,d\mathbf{y}
},
\end{equation}
where $\beta_{\mathrm{src}}=1200$ and $\mathbf{c}$ denotes the coordinate of an injection-production location.

Let $Q_j^{\mathrm{I}}(t)$ and $Q_k^{\mathrm{P}}(t)$ denote the resulting rate at injection and production locations. The source term $g(\mathbf{x},t)$ in the pressure equation~\eqref{eq:thm_pressure} is defined as 
\begin{equation}
    g(\mathbf{x},t)
=
\sum_{j=1}^{12}
Q_j^{\mathrm{I}}(t)
G(\mathbf{x};\mathbf{x}_j^{\mathrm{I}})
-
\sum_{k=1}^{12}
Q_k^{\mathrm{P}}(t)
G(\mathbf{x};\mathbf{x}_k^{\mathrm{P}}).
\end{equation}
The total injection and production rate are balanced every time $\sum_{j=1}^{12}Q_j^{\mathrm{I}}(t)
=
\sum_{k=1}^{12}Q_k^{\mathrm{P}}(t)
=
Q(t).$

The source term $\eta(\mathbf{x},t)$ in the temperature equation \eqref{eq:thm_temperature} is defined by 
$$
\eta(\mathbf{x},t)
=
-\chi_{\mathrm{cold}}
\sum_{j=1}^{12}
Q_j^{\mathrm{I}}(t)
G(\mathbf{x};\mathbf{x}_j^{\mathrm{I}}),
\qquad
\chi_{\mathrm{cold}}=0.85.
$$
This source term represents the cooling caused by fluid injection. Here, $\chi_{\mathrm{cold}}$ controls the strength of the cooling, and the negative
sign indicates the decrease in temperature caused by the
injected cold fluid.

\subsubsection{Numerical and reconstruction setup}
We discretize \(\Omega\) using a uniform \(64\times64\) mesh and use
continuous piecewise-linear finite elements for displacement, pressure, and
temperature. The resulting monolithic system has \(N_h=16900\) degrees of
freedom. We use the time step \(\Delta t=0.01\), giving \(N_t=200\)
solution snapshots over time interval \(\mathcal{J}=(0,2]\). 
We compare POD and SMN for reconstructing the full coupled THM solution. For
SMN, we construct a candidate dictionary containing $s=60$ POD modes and
select $r=15$ modes. For a fair comparison, the standard POD reconstruction
uses the first $r=15$ leading modes.

\subsubsection{Reconstruction results}

The set of mode indices selected by SMN for $r=15$ is $\mathcal{I}_{SMN}
=
\{0,3,8,9,4,6,1,11,5,2,7,10,16,15,14\}$. The selected indices are listed in decreasing order of the magnitude
$|\omega_j|$ of the learned mode-selection weights defined in
Subsection~\ref{subsec:mode_selection}, rather than in increasing numerical order.
Thus, their ordering reflects the importance ranking produced by the SMN
mode-selection procedure. The corresponding full-state and component-wise
reconstruction errors are reported in
Table~\ref{tab:exp1-reconstruction-errors}. 
Table~\ref{tab:exp1-reconstruction-errors} shows that SMN produces smaller
reconstruction errors than POD for the full THM state and every individual
solution component. In particular, the full-state relative error decreases
from $5.8200\times10^{-2}$ for POD to $8.2513\times10^{-5}$ for SMN. 

\begin{table}[!h]
\centering
\renewcommand{\arraystretch}{1.10}
\caption{Experiment 1. Full-state and component-wise relative reconstruction errors for POD and SMN using \(r=15\) modes. The SMN candidate dictionary contains \(s=60\) POD modes.} \label{tab:exp1-reconstruction-errors}
\begin{tabular}{|c|c|c|c|}
\hline
Quantity & Error & \(\mathcal{M}=\mathrm{POD}\)
& \(\mathcal{M}=\mathrm{SMN}\) \\
\hline
Full THM state   & \(E^{\mathcal{M}}\)
& \(5.8200\times10^{-2}\) & \(8.2513\times10^{-5}\) \\
\(u_x\)          & \(E_{u_x}^{\mathcal{M}}\)
& \(1.4850\times10^{-2}\) & \(5.0552\times10^{-6}\) \\
\(u_y\)          & \(E_{u_y}^{\mathcal{M}}\)
& \(3.5665\times10^{-2}\) & \(1.4074\times10^{-5}\) \\
\(p\)            & \(E_p^{\mathcal{M}}\)
& \(5.5365\times10^{-2}\) & \(8.2406\times10^{-5}\) \\
\(\theta\)       & \(E_\theta^{\mathcal{M}}\)
& \(6.5025\times10^{-2}\) & \(8.2794\times10^{-5}\) \\
\(|\mathbf{u}|\) & \(E_{|\mathbf{u}|}^{\mathcal{M}}\)
& \(1.4504\times10^{-2}\) & \(5.0058\times10^{-6}\) \\
\hline
\end{tabular}
\end{table}

We next examine how the reconstruction accuracy changes with the number of
selected modes. We keep the candidate dictionary size fixed at $s=60$ and
consider $r\in\{2,4,6,8,12\}$.
SMN is trained separately for each value of $r$, while the corresponding POD
approximation uses the first $r$ leading modes. The SMN-selected mode-index
sets and the corresponding full-state relative reconstruction errors are
reported in Table~\ref{tab:exp1-multiple-r-errors}. 
Table~\ref{tab:exp1-multiple-r-errors} shows that SMN gives a smaller
full-state reconstruction error than POD for every considered value of $r$.
The improvement becomes especially pronounced as $r$ increases. For example,
when $r=12$, the error decreases from $8.7495\times10^{-2}$ for POD to
$1.2051\times10^{-4}$ for SMN.

\begin{table}[!h]
\centering
\small
\renewcommand{\arraystretch}{1.10}
\setlength{\tabcolsep}{5pt}
\caption{Experiment 1. Selected modes and full-state reconstruction errors for different values of \(r\), with \(s=60\).}\label{tab:exp1-multiple-r-errors}
\begin{tabular}{|c|c|c|c|}
\hline
\(r\) & Selected modes \(I_{\mathrm{SMN}}\)
& \(E^{\mathrm{POD}}\) & \(E^{\mathrm{SMN}}\) \\
\hline
2
& \(\{0,3\}\)
& \(4.1495\times10^{-1}\)
& \(2.9059\times10^{-1}\) \\
\hline
4
& \(\{0,3,4,11\}\)
& \(2.9327\times10^{-1}\)
& \(2.1292\times10^{-1}\) \\
\hline
6
& \(\{0,3,8,11,4,9\}\)
& \(2.0795\times10^{-1}\)
& \(4.8769\times10^{-2}\) \\
\hline
8
& \(\{0,3,8,10,7,9,11,4\}\)
& \(1.5444\times10^{-1}\)
& \(1.6487\times10^{-3}\) \\
\hline
12
& \(\{0,3,8,16,14,4,10,11,7,15,9,2\}\)
& \(8.7495\times10^{-2}\)
& \(1.2051\times10^{-4}\) \\
\hline
\end{tabular}
\end{table}

\begin{figure}[!h] 
\centering 
\begin{subfigure}[t]{\textwidth} 
\centering \includegraphics[width=\textwidth] {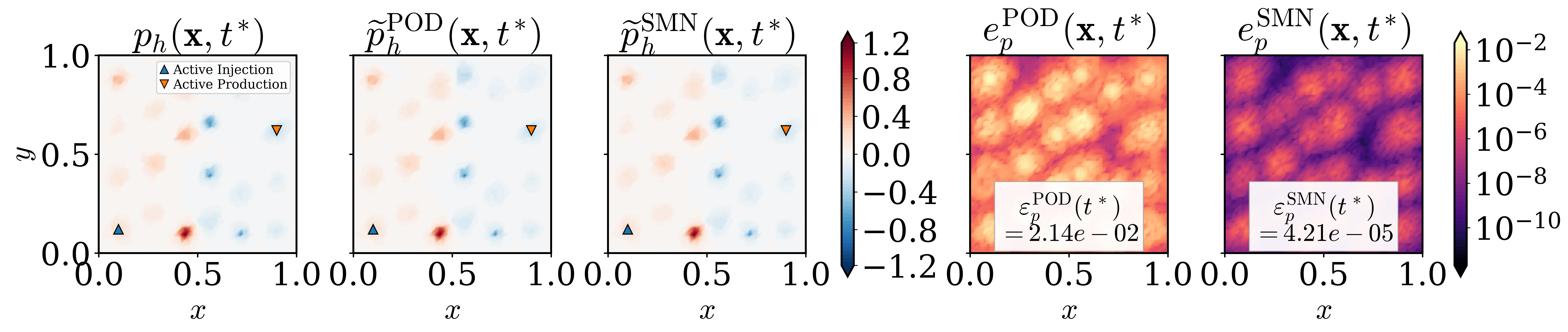} \caption{Pressure reconstruction.} \label{fig:exp1-pressure-reconstruction} \end{subfigure} 
\vspace{0.2em} 
\begin{subfigure}[t]{\textwidth} 
\centering 
\includegraphics[width=\textwidth] {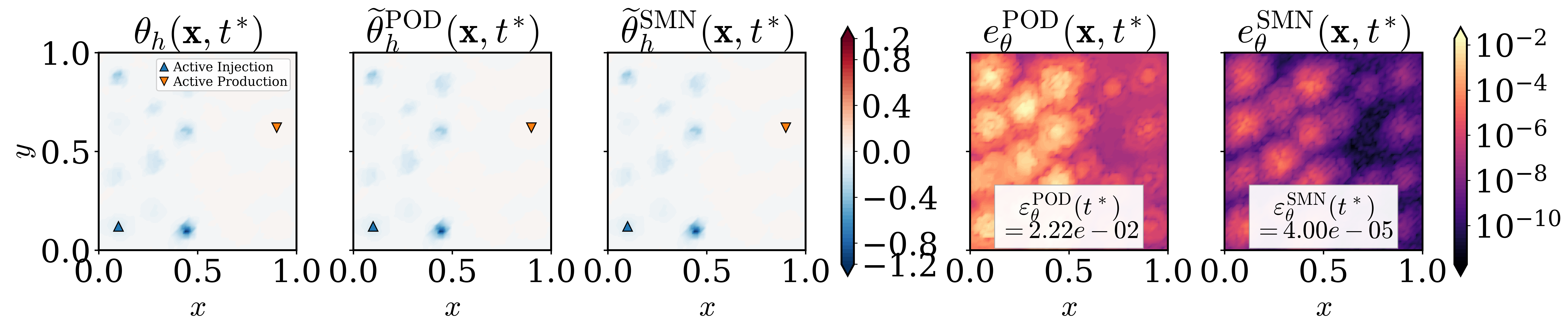} 
\caption{Temperature reconstruction.} \label{fig:exp1-temperature-reconstruction} \end{subfigure} 
\vspace{0.2em} 
\begin{subfigure}[t]{\textwidth} 
\centering \includegraphics[width=\textwidth] {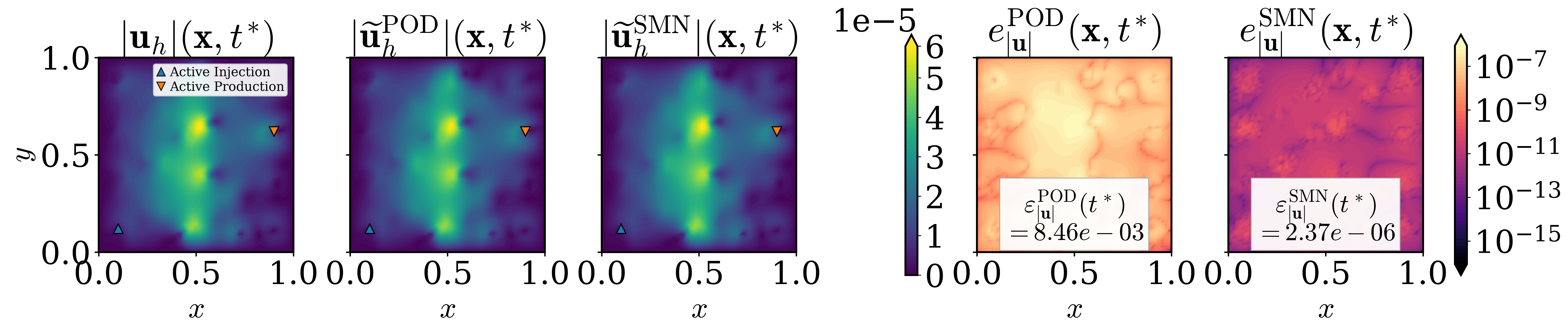} 
\caption{Displacement-magnitude reconstruction.} \label{fig:exp1-displacement-reconstruction} \end{subfigure} 
\caption{Experiment 1. High-fidelity, leading-\(r\) POD, and SMN reconstructions at \(t^*=1.33\), using \(s=60\) candidate modes and \(r=15\) selected modes. The last two panels in each subfigure show the pointwise absolute errors, and the displayed values are the corresponding snapshot relative errors.} \label{fig:exp1-thm-reconstruction} 
\end{figure}

Figure~\ref{fig:exp1-thm-reconstruction} compares the high-fidelity fields
with their POD and SMN reconstructions at the representative time
$t^*=1.33$. For each solution quantity, the first three panels show the
high-fidelity field and its POD and SMN reconstructions, while the last two
panels show the corresponding pointwise absolute errors on a logarithmic
color scale. The values displayed above the error panels are the snapshot
relative errors defined in
equation~\eqref{eq:snapshot-error-measures}. 
The SMN reconstructions closely match the high-fidelity fields and produce
smaller errors than POD. For pressure, temperature, and displacement
magnitude, the SMN snapshot relative errors are
$4.21\times10^{-5}$, $4.00\times10^{-5}$, and $2.37\times10^{-6}$,
respectively. The corresponding POD errors are
$2.14\times10^{-2}$, $2.22\times10^{-2}$, and $8.46\times10^{-3}$.

\subsection{Temporal interpolation using partial training data}
\label{sec:exp1-temporal-interpolation}

In the previous results, all \(N_t=200\) snapshots are used to train and
evaluate SMN. We now examine whether SMN can reconstruct snapshots that are
not included during training. 
We first consider an \(80\%/20\%\) training-test split. We select \(160\)
snapshots at approximately uniformly spaced times over the complete time
interval for training. The remaining \(40\) snapshots are used only for
testing. The POD basis, important modes, and nonlinear decoder are constructed
using only the training snapshots. 
We use $s=60$ candidates POD modes and allow SMN to select $r=15$ modes. For the training case, the selected mode indices are $\mathcal{I}_{\mathrm{SMN}}
=
\{0,5,4,7,6,3,1,2,11,8,14,24,21,35,31\}.$ During testing, each unseen snapshot is projected onto these modes. The
resulting \(15\) coefficients are passed through the trained nonlinear
decoder to reconstruct the complete THM solution. The selected modes and
decoder parameters remain fixed during testing. The training and test errors
are reported in Table~\ref{tab:exp1-uniform-train-test}. For the unseen test snapshots, the SMN error is
\(1.0625\times10^{-4}\), compared with
\(3.4261\times10^{-2}\) for POD. 

\begin{table}[!h]
\centering
\renewcommand{\arraystretch}{1.10}
\caption{Experiment 1. Full-state reconstruction errors for the uniformly distributed \(80\%/20\%\) training-test split. We use \(s=60\) candidate POD modes and \(r=15\) selected modes. Only the training snapshots are used to construct the POD basis and train SparseModesNet.}
\label{tab:exp1-uniform-train-test}
\begin{tabular}{|c|c|c|}
\hline
Data set & \(E^{\mathrm{POD}}\) & \(E^{\mathrm{SMN}}\) \\
\hline
Training snapshots
& \(3.4292\times10^{-2}\)
& \(3.2541\times10^{-5}\) \\
Test snapshots
& \(3.4261\times10^{-2}\)
& \(1.0625\times10^{-4}\) \\
\hline
\end{tabular}
\end{table}

We next examine the effect of using less training data. We repeat the
experiment using \(40\%\) and \(60\%\) of the snapshots for training,
corresponding to \(80\) and \(120\) training snapshots, respectively. All
remaining snapshots are used for testing. For each case, the POD basis,
selected modes, and nonlinear decoder are recomputed using only the
corresponding training snapshots. The values \(s=60\) and \(r=15\) remain
fixed. 
Table~\ref{tab:exp1-uniform-fraction-study} shows that SMN test error decreases from \(7.5116\times10^{-3}\) with \(40\%\)
training data to \(3.8997\times10^{-3}\) with \(60\%\) training data. For the
\(80\%\) case, it further decreases to \(1.0625\times10^{-4}\). SMN also
gives a smaller test error than POD for every training fraction considered.
These results show that SMN can interpolate the THM solution at unseen times
and that its accuracy improves when more training snapshots are available.

\begin{table}[!h]
\centering
\renewcommand{\arraystretch}{1.10}
\caption{Experiment 1. Test reconstruction errors obtained using \(40\%\) and \(60\%\) of the snapshots for training. The values \(s=60\) and \(r=15\) are fixed in both cases.}
\label{tab:exp1-uniform-fraction-study}
\begin{tabular}{|c|c|c|c|c|}
\hline
Training fraction
& Training snapshots
& Test snapshots
& \(E^{\mathrm{POD}}\)
& \(E^{\mathrm{SMN}}\) \\
\hline
\(40\%\)
& \(80\)
& \(120\)
& \(3.4381\times10^{-2}\)
& \(7.5116\times10^{-3}\) \\
\(60\%\)
& \(120\)
& \(80\)
& \(3.4312\times10^{-2}\)
& \(3.8997\times10^{-3}\) \\
\hline
\end{tabular}
\end{table}

\subsection{Experiment 2: THM system with fracture-network coefficients and multiple well locations}
\label{subsec:exp2_ghs}

\subsubsection{Experiment setup and fracture-controlled circulation}

In the second experiment, we consider THM problem with a prescribed fracture network and several fixed injection and production locations. Unlike Experiment \ref{subsec:exp1_thm_heterogeneous}, in which the permeability and thermal conductivity coefficient values vary randomly from one mesh element to another, the heterogeneity in this experiment follows the geometry of the connected fracture network. The purpose of this experiment is to test whether SMN can accurately reconstruct the fracture-controlled THM solution and preserve selected geothermal quantities of interest.

\begin{figure}[!h] 
\centering 
\begin{subfigure}[t]{0.48\textwidth} 
\centering \includegraphics[width=\textwidth] {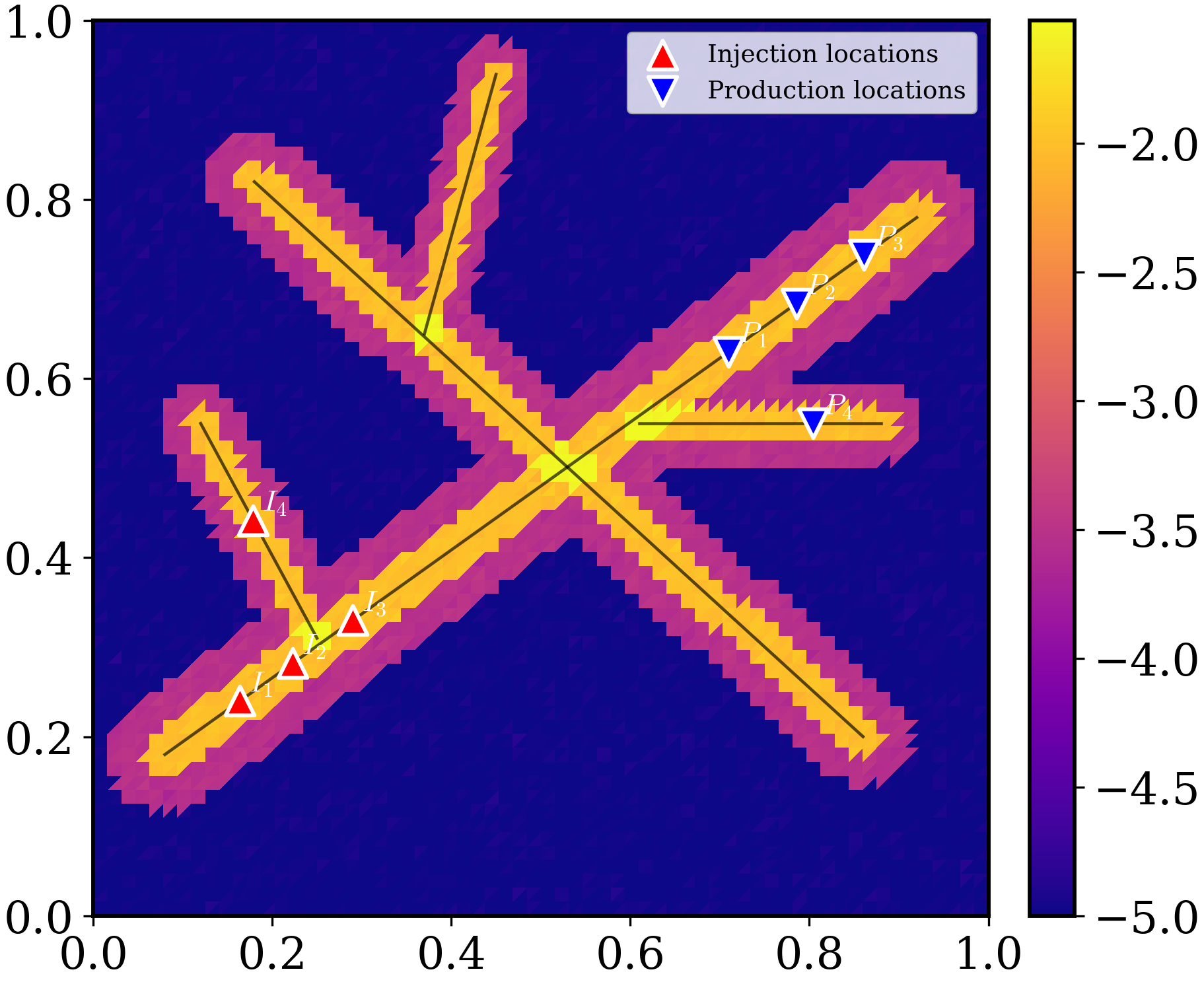} \caption{Permeability \(\log_{10}\kappa(\mathbf{x})\).} \label{fig:exp2-logK-multifeed} \end{subfigure} \hfill \begin{subfigure}[t]{0.48\textwidth} \centering \includegraphics[width=\textwidth] {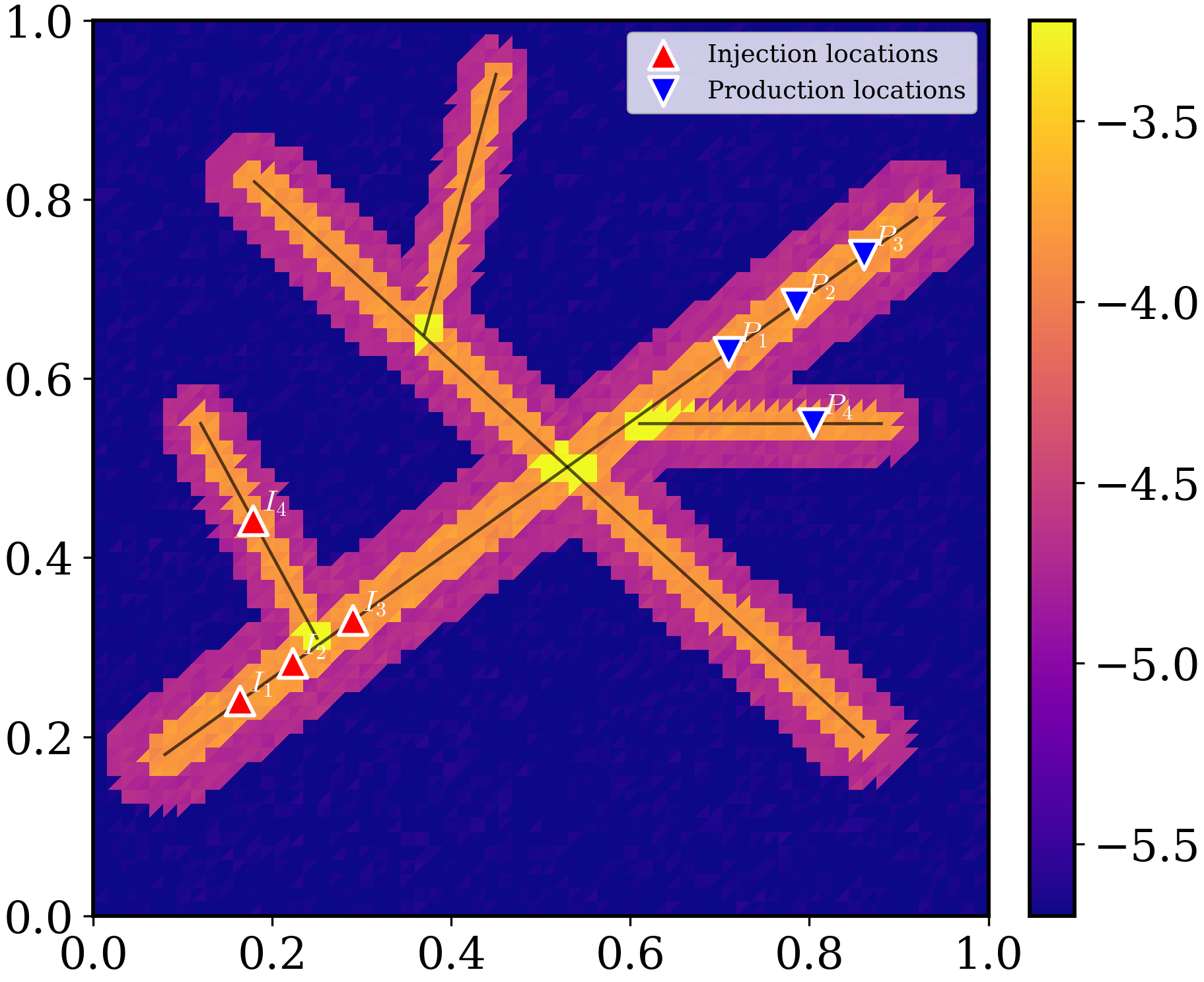} \caption{Thermal conductivity \(\log_{10}d_T(\mathbf{x})\).} 
\label{fig:exp2-logD-multifeed} 
\end{subfigure} 
\caption{Experiment 2. Fracture-network coefficient fields and fixed injection and production locations. The injection locations \(I_1,\ldots,I_4\) and production locations \(P_1,\ldots,P_4\) are marked by red upward and blue downward triangles, respectively.} \label{fig:exp2-fracture-kd-fields} 
\end{figure}

The permeability and thermal conductivity tensors are isotropic and are written as $\mathbf{K}(\mathbf{x})=\kappa(\mathbf{x})\mathbf{I}_2$ and 
$\mathbf{D}(\mathbf{x})=d_T(\mathbf{x})\mathbf{I}_2$, where $\mathbf{I}_2$ is the $2 \times 2$ identity matrix. Both coefficients have larger values inside the fractures than in the surrounding matrix. The fractures therefore provide preferential pathways for fluid flow and heat transfer and influence the pressure and temperature fields. The coefficient fields and the fixed injection and production locations shown in Figure~\ref{fig:exp2-fracture-kd-fields}. All locations lie on the fracture network and remain fixed throughout the simulation.

We use the same domain, THM parameters, initial and boundary conditions, spatial discretization as in Experiment \ref{subsec:exp1_thm_heterogeneous}. The problem is solved over the time interval $\mathcal{J}=(0,4]$ using the backward Euler method with $\Delta t=0.01$. The monolithic THM state has $N_h=16900$ degrees of freedom, and $N_t=400$ snapshots are collected. Therefore, $\mathbf{Z}\in\mathbb{R}^{16900\times400}$.

To represent geothermal circulation, cold fluid is injected at four fixed injection locations and fluid is removed at four fixed production locations. The connected fractures provide the main pathways between these locations. Similar injection and production
systems in fractured geothermal reservoirs have been studied in
\cite{lei2023injectionproduction}. In this experiment, the common total injection and production rate change with time. This rate is also redistributed among the fixed locations. The redistribution of flow among fixed geothermal wells has been considered
in~\cite{akin2014reinjection,duplyakin2022geothermal}. Time-dependent and
cyclic flow rates have also been used in studies of fractured geothermal
reservoirs~\cite{nakao2003cyclic,fokker2026harmonic,hofmann2018cyclic}.

The circulation schedule has two time-dependent components. First, the total
rate $q_{\mathrm{tot}}(t)$ determines the overall amount of fluid injected into
and removed from the reservoir at time $t$. Second, the weights
$w_j^{\mathrm I}(t)$ and $w_k^{\mathrm P}(t)$ determine how this total rate is
distributed among the four injection and four production locations,
respectively. Thus, $q_{\mathrm{tot}}(t)$ determines how much fluid is
circulated, whereas the weights determine where the fluid is injected and
produced.

Let $\mathbf{x}_j^{\mathrm I}$, $j=1,\ldots,4$, and $\mathbf{x}_k^{\mathrm P}$, $k=1,\ldots,4$, denote the fixed injection and production locations, respectively. The common total injection and production rate is defined by 
$$
q_{\mathrm{tot}}(t) = q_0\left[ 1+a_1\sin\left(\frac{2\pi t}{T_q}+\phi_1\right) +a_2\sin\left(\frac{4\pi t}{T_q}+\phi_2\right) \right],
$$
where \(q_0=3\), \(a_1=0.25\), \(a_2=0.10\), \(T_q=1.35\), \(\phi_1=0\), and \(\phi_2=0.35\pi\). The total rate is distributed among the four injection and four production locations using non-negative time-dependent weights. To define these weights, we first introduce 
$$
\widehat{w}(t;\psi)
=
b_w +
a_w
\left[
1+\sin\left(
\frac{2\pi t}{T_w}+\psi
\right)
\right],
$$
where, $b_w=0.12, a_w=0.10,T_w=4$ and $\psi$ is a location-dependent phase shift. The phase shift used for the injection and production location are $\left(\psi_1^{\mathrm I},\ldots,\psi_4^{\mathrm I}\right)
=
(0,0.45\pi,0.95\pi,1.35\pi),$ and $\left(\psi_1^{\mathrm P},\ldots,\psi_4^{\mathrm P}\right)
=
(0.25\pi,0.75\pi,1.15\pi,1.65\pi),$ respectively. The normalized injection and production weights are then defined by
$$
w_j^{\mathrm I}(t)
=
\frac{
\widehat{w}(t;\psi_j^{\mathrm I})
}{
\displaystyle
\sum_{\ell=1}^{4}
\widehat{w}(t;\psi_{\ell}^{\mathrm I})
},
\qquad j=1,\ldots,4,
$$
and 
$$
w_k^{\mathrm P}(t)
=
\frac{
\widehat{w}(t;\psi_k^{\mathrm P})
}{
\displaystyle
\sum_{\ell=1}^{4}
\widehat{w}(t;\psi_{\ell}^{\mathrm P})
},
\qquad k=1,\ldots,4.
$$

Because, $\widehat{w}(t;\psi)\geq b_w>0$, these weights are well-defined and non-negative. Moreover,   these weights satisfy $\sum_{j=1}^4 w_j^{\mathrm I}(t)=1$, and $\sum_{k=1}^4 w_k^{\mathrm P}(t)=1$.

The rates assigned to each location are 
\begin{equation}\label{eq:exp3-in-prod-rate}
q_j^{\mathrm I}(t) = q_{\mathrm{tot}}(t)w_j^{\mathrm I}(t) \qquad \mbox{and} \qquad q_k^{\mathrm P}(t) = q_{\mathrm{tot}}(t)w_k^{\mathrm P}(t). 
\end{equation}
Therefore, the total injection and production rate remain balanced while their distributions among the fixed locations change with time. Each location is represented by the localized Gaussian function defined in \eqref{eq:gaussian_fun}, where $\beta_{\mathrm{src}}=1200$ and $\mathbf{c}$ denotes the coordinates of the corresponding injection or production location.

We then set $
G_j^{\mathrm I}(\mathbf{x})
=
G(\mathbf{x};\mathbf{x}_j^{\mathrm I})$, and 
$G_k^{\mathrm P}(\mathbf{x})
=
G(\mathbf{x};\mathbf{x}_k^{\mathrm P}).
$ The pressure and temperature source terms are 
$$
\begin{aligned} 
g(\mathbf{x},t) &= \sum_{j=1}^4 q_j^{\mathrm I}(t)G_j^{\mathrm I}(\mathbf{x}) - \sum_{k=1}^4 q_k^{\mathrm P}(t)G_k^{\mathrm P}(\mathbf{x}),\\ \eta(\mathbf{x},t) &= -\chi_{\mathrm{cold}} \sum_{j=1}^4 q_j^{\mathrm I}(t)G_j^{\mathrm I}(\mathbf{x}), \qquad \chi_{\mathrm{cold}}=0.85. \end{aligned} $$
Thus, $g$ represents injection and production, while $\eta$ represent cooling near the injection locations. As \(q_{\mathrm{tot}}(t)\) and the weights vary, the strength and relative contribution of the different locations change with time. This produces pressure and temperature patterns that evolve through different parts of the fracture network.

\subsubsection{Geothermal quantities of interest}
In addition to reconstruction errors, we evaluate two geothermal quantities of interest. The first measures the average temperature change near the production locations. Production temperature is commonly used to
evaluate the thermal performance of a geothermal system ~\cite{lei2023injectionproduction}. Because wells are not explicitly resolved in our model, we compute a localized average around each production location. 

Let \(\{\mathbf{x}_i\}_{i=1}^{N_v}\) denote the mesh vertices. For a scalar
finite element function \(v_h\), we define the localized average near
\(P_k\) by
\[
\mathcal{A}_k^{\mathrm P}(v_h)
=
\frac{
\displaystyle\sum_{i=1}^{N_v}
e^{-\gamma\|\mathbf{x}_i-\mathbf{x}_k^{\mathrm P}\|_2^2}
v_h(\mathbf{x}_i)
}{
\displaystyle\sum_{i=1}^{N_v}
e^{-\gamma\|\mathbf{x}_i-\mathbf{x}_k^{\mathrm P}\|_2^2}
},
\qquad
\gamma=\frac{\beta_{\mathrm{src}}}{2}=600.
\]
We combine the four localized averages using the production weights: 
$$
\overline{\theta}_{\mathrm P}(t^n) = \sum_{k=1}^4 w_k^{\mathrm P}(t^n) \mathcal{A}_k^{\mathrm P} \bigl(\theta_h(\cdot,t^n)\bigr). $$
The first quantity of interest is then defined as 
$\mathcal{Q}_{\theta}(t^n) = \overline{\theta}_{\mathrm P}(t^n)-
\overline{\theta}^{\mathrm P}(0).$ The quantity represents the temperature change near the production locations relative to initial temperature. Negative values indicate cooling, positive values indicate heating, and values close to zero indicate little temperature change.

The second quantity of interest measures the pressure difference between injection and production locations. In geothermal studies, the injection-production pressure difference divided by the circulation rate used to describe the hydraulic response of the reservoir ~\cite{zeng2021welllayout}. In our model, we use corresponding pressure difference divided by the prescribed total source amplitude. 

For a scalar finite element function $v_h$, we define its localized average near the injection location $I_j$ by 
$$\mathcal{A}_j^{\mathrm I}(v_h) = \frac{ \displaystyle\sum_{i=1}^{N_v} e^{-\gamma\|\mathbf{x}_i-\mathbf{x}_j^{\mathrm I}\|_2^2} v_h(\mathbf{x}_i) }{ \displaystyle\sum_{i=1}^{N_v} e^{-\gamma\|\mathbf{x}_i-\mathbf{x}_j^{\mathrm I}\|_2^2} }, \qquad \gamma=600. $$
Using the injection and production weights, we define the corresponding pressure averages by 
$$\overline{p}_{\mathrm I}(t^n) = \sum_{j=1}^4 w_j^{\mathrm I}(t^n) \mathcal{A}_j^{\mathrm I} \bigl(p_h(\cdot,t^n)\bigr), \qquad \overline{p}_{\mathrm P}(t^n) = \sum_{k=1}^4 w_k^{\mathrm P}(t^n) \mathcal{A}_k^{\mathrm P} \bigl(p_h(\cdot,t^n)\bigr). $$ 
The second quantity of interest is then defined by 
\begin{equation}
\mathcal{Q}_{p}(t^n) = \frac{ \overline{p}_{\mathrm I}(t^n) - \overline{p}_{\mathrm P}(t^n) }{ q_{\mathrm{tot}}(t^n) }. 
\end{equation}
Thus, $\mathcal{Q}_{p}(t^n)$ measures the pressure difference between the
injection and production locations per unit total source amplitude. 

\subsubsection{Reconstruction results}

We apply SMN to the full THM snapshot matrix using $s=100$ candidate POD modes and $r=10$ selected modes. The selected modes indices are 
$\mathcal{I}_{SMN} = \{0,1,5,4,3,2,6,8,16,19\}.$ 
The corresponding full-state and component-wise reconstruction errors are reported in Table~\ref{tab:exp2-reconstruction-errors}, highlighting a significant difference between POD and SMN where SMN outperforms classical POD in the reconstruction of the complex fields at stake. 

\begin{table}[!h]
\centering
\small
\renewcommand{\arraystretch}{1.10}
\caption{Experiment 2. Full-state and component-wise reconstruction errors for \(s=100\) and \(r=10\).}
\label{tab:exp2-reconstruction-errors}
\begin{tabular}{|c|c|c|}
\hline
State or component
& POD error
& SMN error \\
\hline
Full THM state
& \(E^{\mathrm{POD}}=1.3670\times10^{-2}\)
& \(E^{\mathrm{SMN}}=1.0597\times10^{-8}\) \\
\(u_x\)
& \(E_{u_x}^{\mathrm{POD}}=4.0885\times10^{-3}\)
& \(E_{u_x}^{\mathrm{SMN}}=1.8977\times10^{-8}\) \\
\(u_y\)
& \(E_{u_y}^{\mathrm{POD}}=9.8101\times10^{-3}\)
& \(E_{u_y}^{\mathrm{SMN}}=4.8552\times10^{-8}\) \\
\(p\)
& \(E_p^{\mathrm{POD}}=3.4201\times10^{-2}\)
& \(E_p^{\mathrm{SMN}}=3.9856\times10^{-7}\) \\
\(\theta\)
& \(E_\theta^{\mathrm{POD}}=3.4345\times10^{-2}\)
& \(E_\theta^{\mathrm{SMN}}=3.9308\times10^{-7}\) \\
\(|\mathbf{u}|\)
& \(E_{|\mathbf{u}|}^{\mathrm{POD}}=3.3019\times10^{-3}\)
& \(E_{|\mathbf{u}|}^{\mathrm{SMN}}=1.3926\times10^{-8}\) \\
\hline
\end{tabular}
\end{table}

We next examine whether the reconstructed pressure and temperature field preserve two quantities of interest. For $\ell\in\{\theta,p\}$, let $\mathbf{Q}_{\ell} = \left[ \mathcal{Q}_{\ell}(t^1),\ldots, \mathcal{Q}_{\ell}(t^{N_t}) \right]$, and let $\widetilde{\mathbf{Q}}_{\ell}^{\mathcal M}$ denotes the corresponding values computed from the reconstructed method $\mathcal M\in\{\mathrm{POD},\mathrm{SMN}\}$. We measure the pointwise absolute error and the relative error over the complete time interval by
\begin{equation}
    e_{\mathcal{Q}_{\ell}}^{\mathcal M}(t^n) = \left| \mathcal{Q}_{\ell}(t^n) - \widetilde{\mathcal{Q}}_{\ell}^{\mathcal M}(t^n) \right| \qquad \mbox{and} \qquad E_{\mathcal{Q}_{\ell}}^{\mathcal M} = \frac{ \left\| \mathbf{Q}_{\ell} - \widetilde{\mathbf{Q}}_{\ell}^{\mathcal M} \right\|_2 }{ \left\|\mathbf{Q}_{\ell}\right\|_2 }. 
\end{equation}

\begin{figure}[!h]
\centering
\begin{subfigure}[t]{0.47\textwidth}
\centering
\includegraphics[width=\textwidth]
{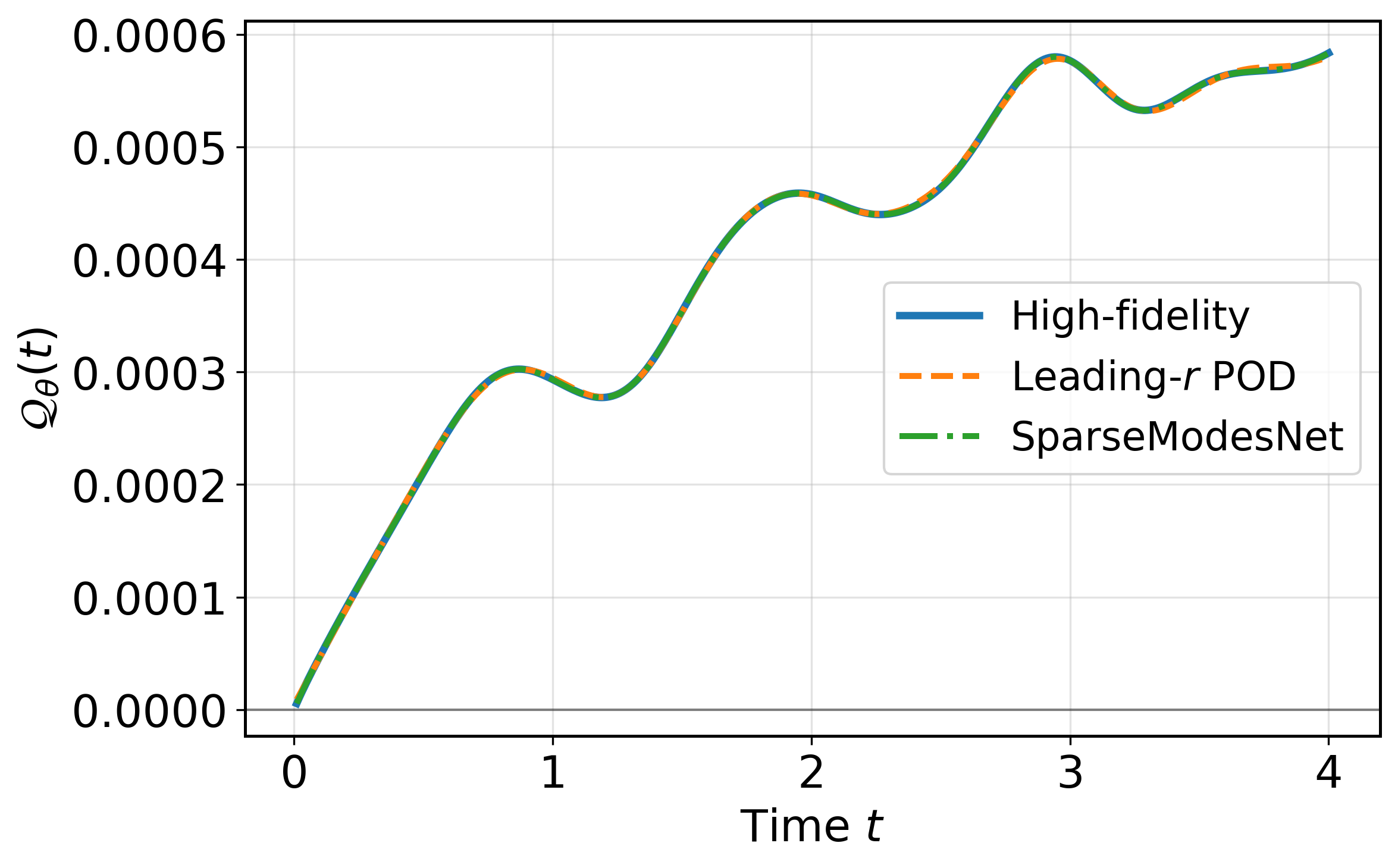}
\caption{Time history of \(\mathcal{Q}_{\theta}\).}
\label{fig:exp2-qoi-temperature-curve}
\end{subfigure}
\hfill
\begin{subfigure}[t]{0.47\textwidth}
\centering
\includegraphics[width=\textwidth]
{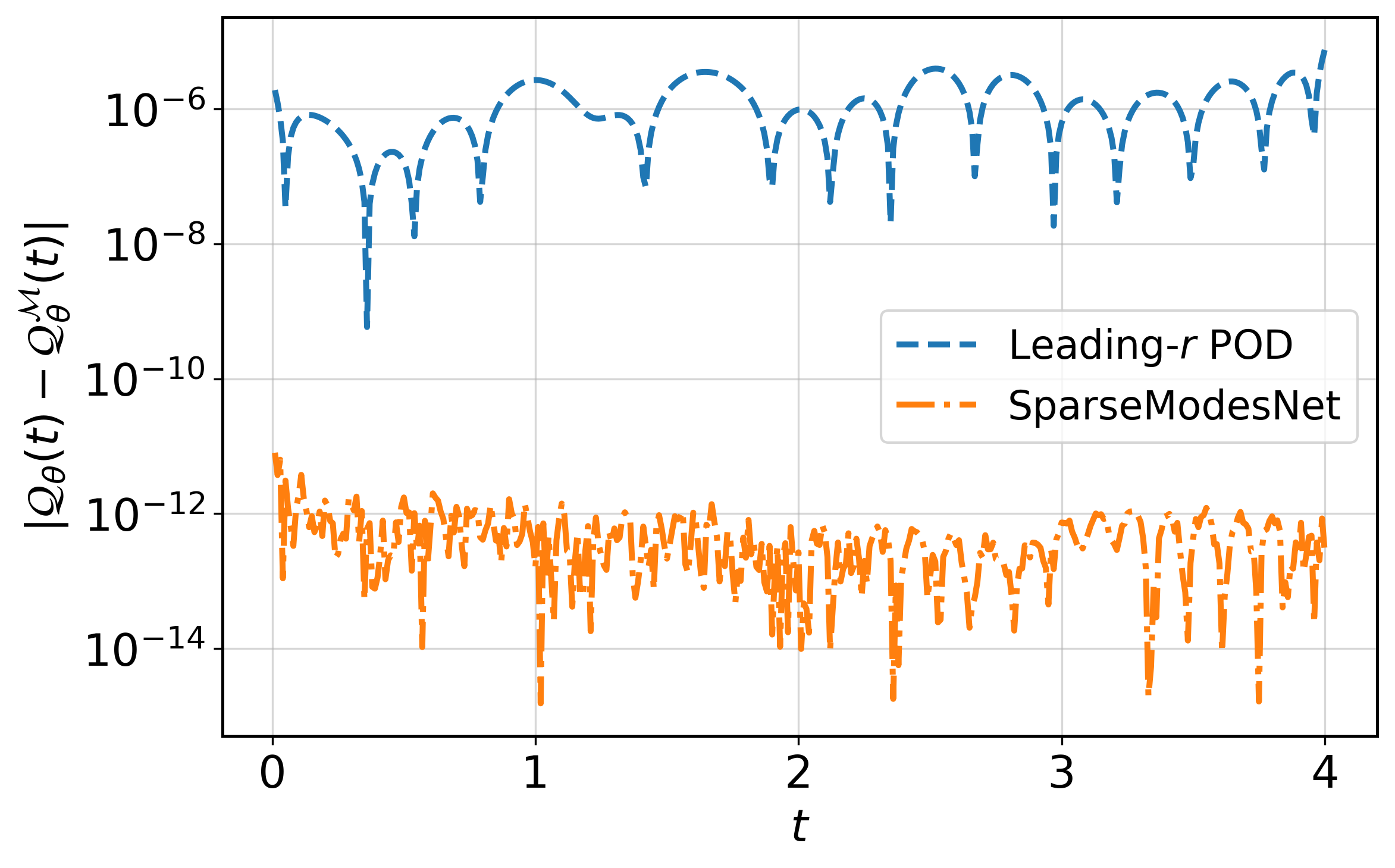}
\caption{Pointwise absolute error.}
\label{fig:exp2-qoi-temperature-error}
\end{subfigure}
\caption{Experiment 2. Reconstruction of the production-side temperature quantity
\(\mathcal{Q}_{\theta}\) using the leading-\(r\) POD and SparseModesNet,
with \(s=100\) and \(r=10\).}
\label{fig:exp2-qoi-temperature}
\end{figure}

Figure~\ref{fig:exp2-qoi-temperature} compares the time histories of $\mathcal{Q}_{\theta}$. The high-fidelity, POD, and SMN curves are visually indistinguishable. However, the absolute-error curves show that SMN is much more accurate than POD. The relative errors are $E_{\mathcal{Q}_{\theta}}^{\mathrm{POD}} = 8.1881\times10^{-3}$, and  $E_{\mathcal{Q}_{\theta}}^{\mathrm{SMN}} = 3.9466\times10^{-9}.$ 
Thus, both methods reproduce the overall time dependence of \(\mathcal{Q}_{\theta}\), but SMN gives a much smaller error.

\begin{figure}[!h]
\centering
\begin{subfigure}[t]{0.47\textwidth}
\centering
\includegraphics[width=\textwidth]
{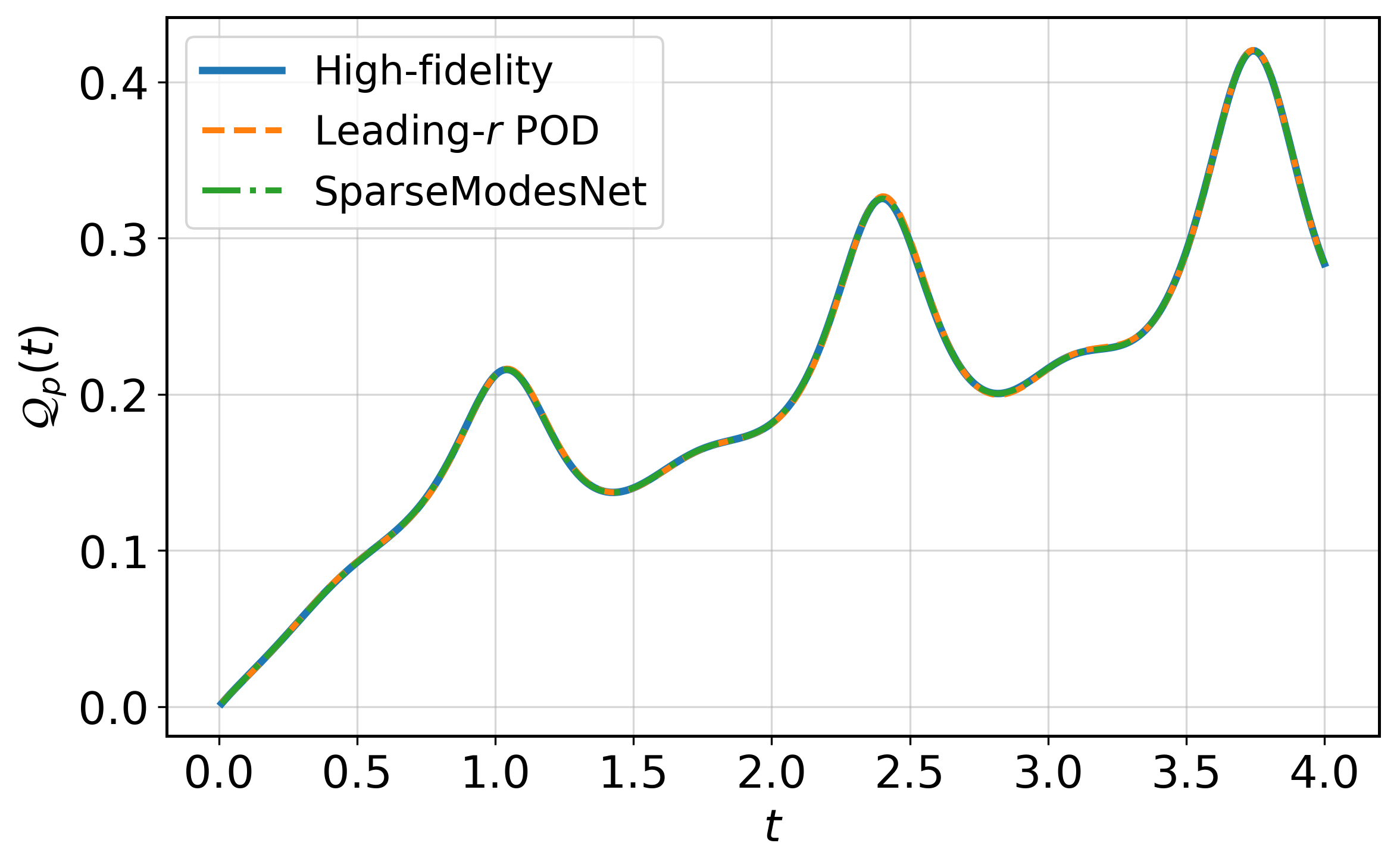}
\caption{Time history of \(\mathcal{Q}_{p}\).}
\label{fig:exp2-qoi-pressure-curve}
\end{subfigure}
\hfill
\begin{subfigure}[t]{0.47\textwidth}
\centering
\includegraphics[width=\textwidth]
{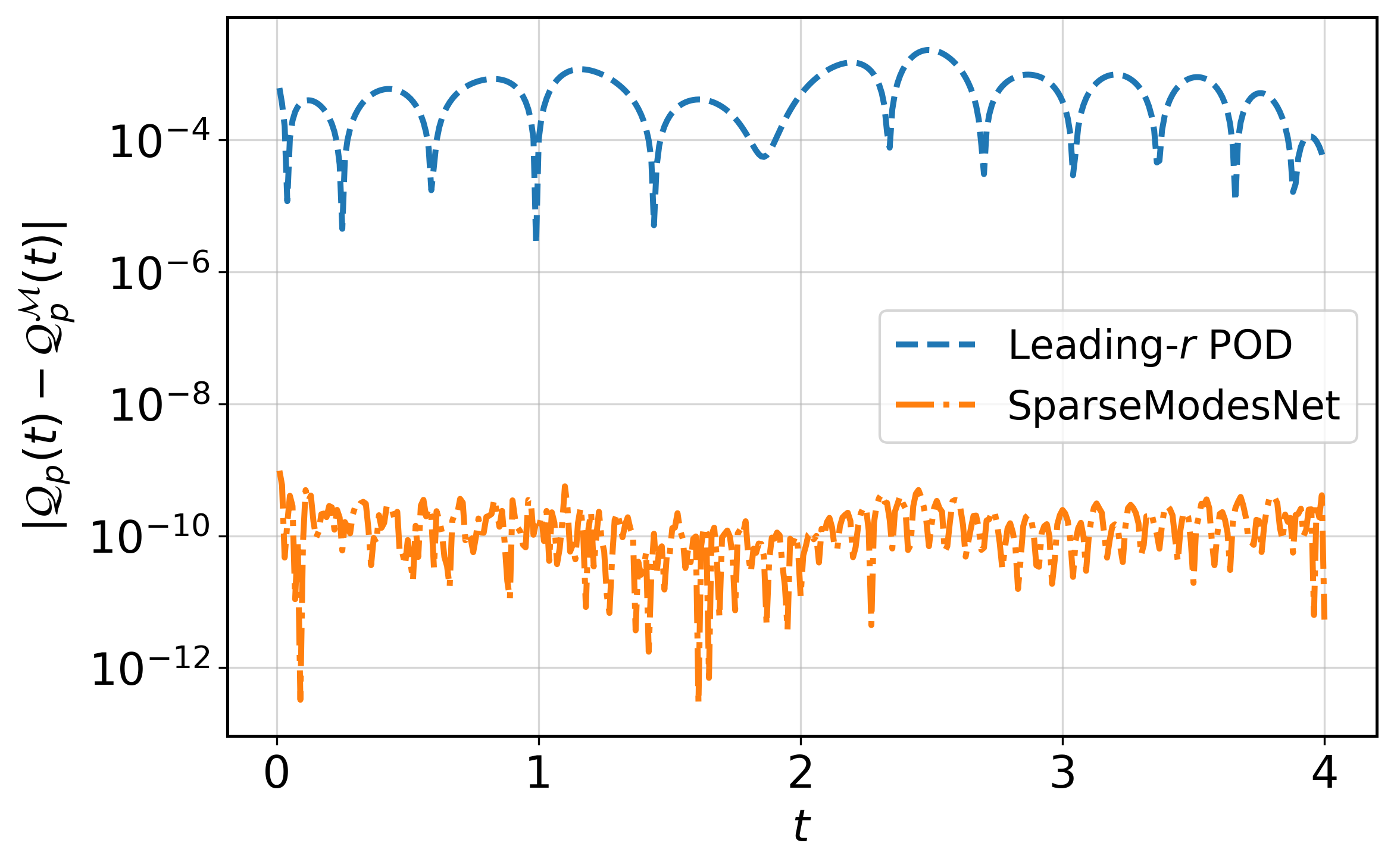}
\caption{Pointwise absolute error.}
\label{fig:exp2-qoi-pressure-error}
\end{subfigure}
\caption{Experiment 2. Reconstruction of the source-normalized injection--production
pressure difference \(\mathcal{Q}_{p}\) using the leading-\(r\) POD and
SparseModesNet, with \(s=100\) and \(r=10\).}
\label{fig:exp2-qoi-pressure}
\end{figure}

Figure~\ref{fig:exp2-qoi-pressure} presents the corresponding results for \(\mathcal{Q}_{p}\). The three time histories again appear nearly identical, but the absolute-error curves clearly show the difference between the two methods. The relative errors are $ E_{\mathcal{Q}_{p}}^{\mathrm{POD}} = 3.45\times10^{-3}$, and $E_{\mathcal{Q}_{p}}^{\mathrm{SMN}} = 8.71\times10^{-10}.$ Therefore, SMN reconstructs the time-dependent pressure quantity more accurately than the leading-r POD approximation.  
Overall, both POD and SMN reproduce the visible behavior of the two quantities of interest. However, the relative errors obtained with SMN are approximately six orders of magnitude smaller than those obtained with POD.

\subsection{Experiment 3: Three-dimensional THM system with Utah FORGE
fracture geometry}
\label{subsec:exp3_forge}
\subsubsection{Utah FORGE fracture geometry and coefficient fields}

In the third experiment, we extend the two-dimensional fracture-network problem considered in Experiment \ref{subsec:exp2_ghs} to three dimensions. The fracture geometry is based on the Utah FORGE 2025 v1 discrete fracture network data set~\cite{finnila2025utahforge_dfn}. The purpose of the this experiment to test whether SMN can reconstruct a three-dimensional THM solution  and preserve pressure and temperature based quantities of interest.

\begin{figure}[!h]
\centering
\includegraphics[width=0.60\textwidth]
{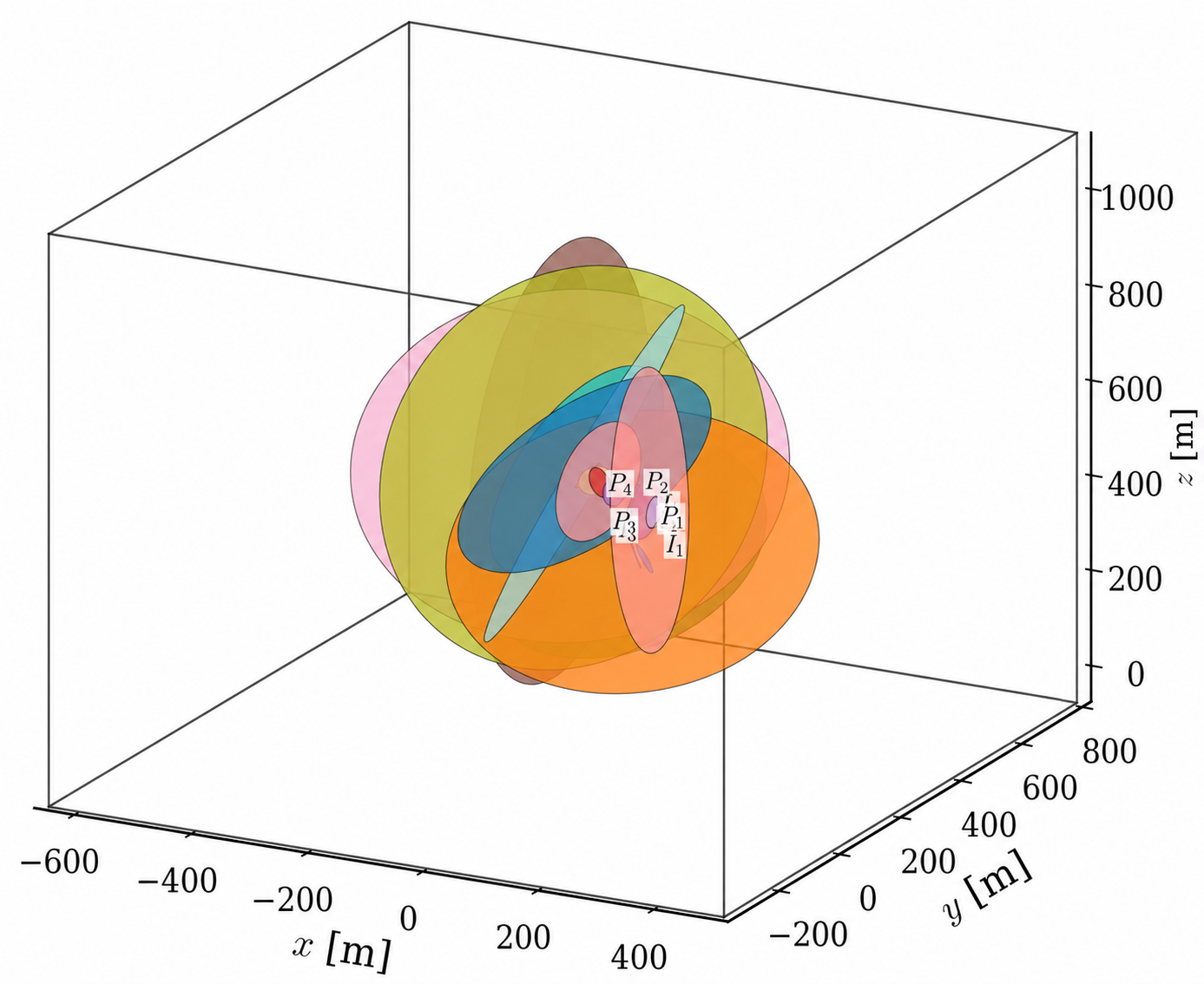}
\caption{Experiment 3. Utah FORGE Stage~8 Case~A+B fracture geometry. Each fracture is represented here by a planar circular disk
constructed from its reported center, orientation, and radius. $I_1,\ldots,I_4$, and 
$P_1,\ldots,P_4$ represented the injection and production locations, respectively. The rectangular box represents the computational domain.
The disk colors distinguish the individual fractures and do not represent
coefficient values.}
\label{fig:exp3-stage8-geometry}
\end{figure}

We use the Stage~8 Case~A+B subset of the dataset, which contains $N_f=31$ fractures: $15$ Case~A fractures and $16$ Case~B fractures. For each fracture $\ell=1,\ldots,N_f$, the data provide its center $\mathbf{c}_\ell$, radius $R_\ell$, and orientation. The reported pole trend and plunge are used to construct the unit normal vector $\mathbf{n}_\ell$. The Utah FORGE reference DFN assumes that the fractures are approximately circular ~\cite{finnila2019building}. Planar circular disks are also commonly used in DFN models~\cite{gomez2023nonparametric}. We therefore represent fracture $\ell$ by the planar disk
\begin{equation}
    F_\ell = \left\{ \mathbf{x}\in\mathbb{R}^3: (\mathbf{x}-\mathbf{c}_\ell)\cdot\mathbf{n}_\ell=0, \quad \|\mathbf{x}-\mathbf{c}_\ell\|_2\leq R_\ell \right\}. 
\end{equation}
This representation preserves the reported centers, size and orientation of the fracture, but it does not describe its exact physical boundary. Figure~\ref{fig:exp3-stage8-geometry} shows the selected fracture geometry inside computational domain. The four injection location $I_1,\ldots I_4$ and four production locations $P_1,\ldots P_4$, which are used later to define the source terms, are also shown. 

We next use the fracture geometry to construct the permeability and thermal conductivity fields. Because each planar disk $F_\ell$ has zero thickness, it
cannot be represented directly by the three-dimensional volumetric mesh. We
therefore assign each fracture a numerical thickness
$\delta_f=10\,\mathrm{m}$ and define the fracture-supported region by
\begin{equation}
    F_\ell^{\delta_f}
=
\left\{
\mathbf{x}\in\Omega:
\left|(\mathbf{x}-\mathbf{c}_\ell)\cdot\mathbf{n}_\ell\right|
\leq\frac{\delta_f}{2},
\quad
\left\|
(\mathbf{I}_3-\mathbf{n}_\ell\mathbf{n}_\ell^\top)
(\mathbf{x}-\mathbf{c}_\ell)
\right\|_2
\leq R_\ell
\right\}.
\end{equation}
where $\mathbf{I}_3$ is the $3\times3$ identity matrix. The first condition selects points that are close to the fracture plane. The second condition ensures that these points remain within the radius of the fracture. The thickness \(\delta_f\) is used only to represent the fracture on the computational mesh.

\begin{figure}[!h]
\centering
\includegraphics[width=0.98\textwidth]
{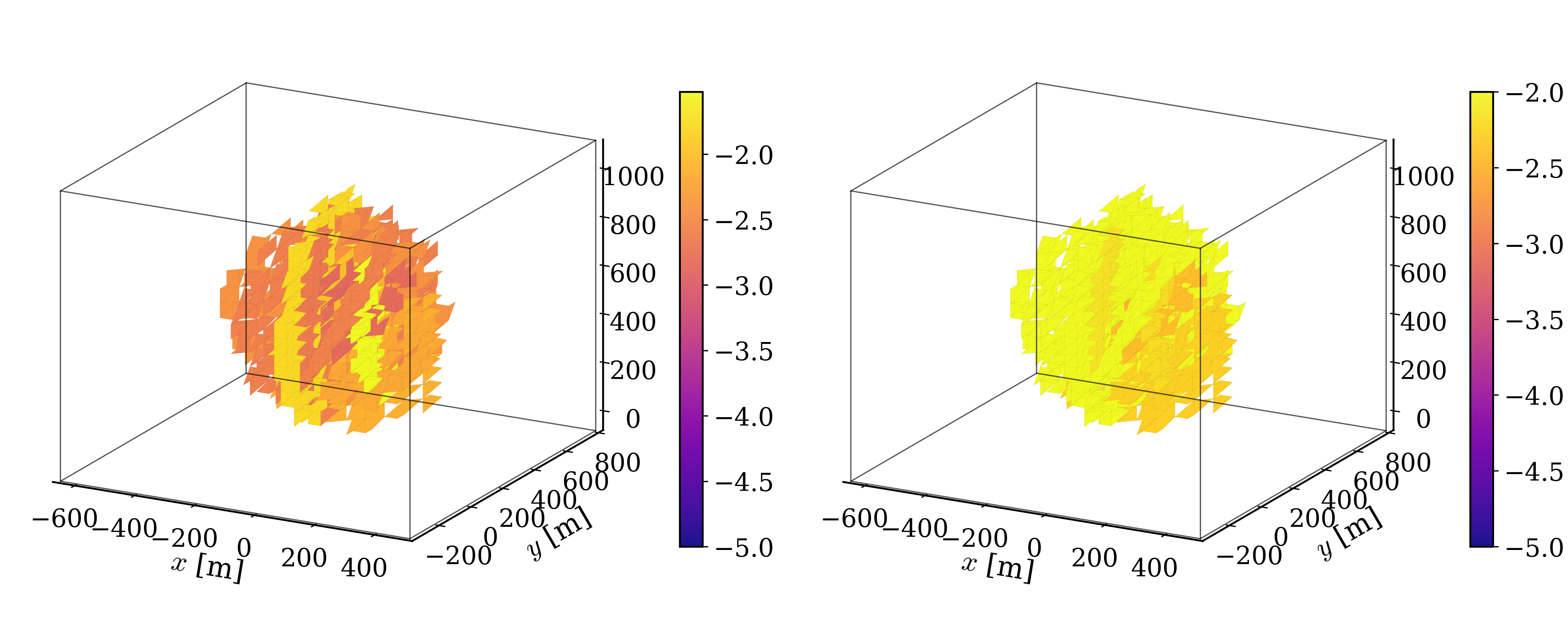}
\caption{Experiment 3. Three-dimensional view of the cellwise-constant coefficient fields:
(a) permeability \(\log_{10}\kappa\) and (b) thermal conductivity
\(\log_{10}d_T\). Only the fracture-supported elements are displayed.}
\label{fig:exp3-computational-coefficients-3d}
\end{figure}

\begin{figure}[!h]
\centering
\includegraphics[width=0.95\textwidth]
{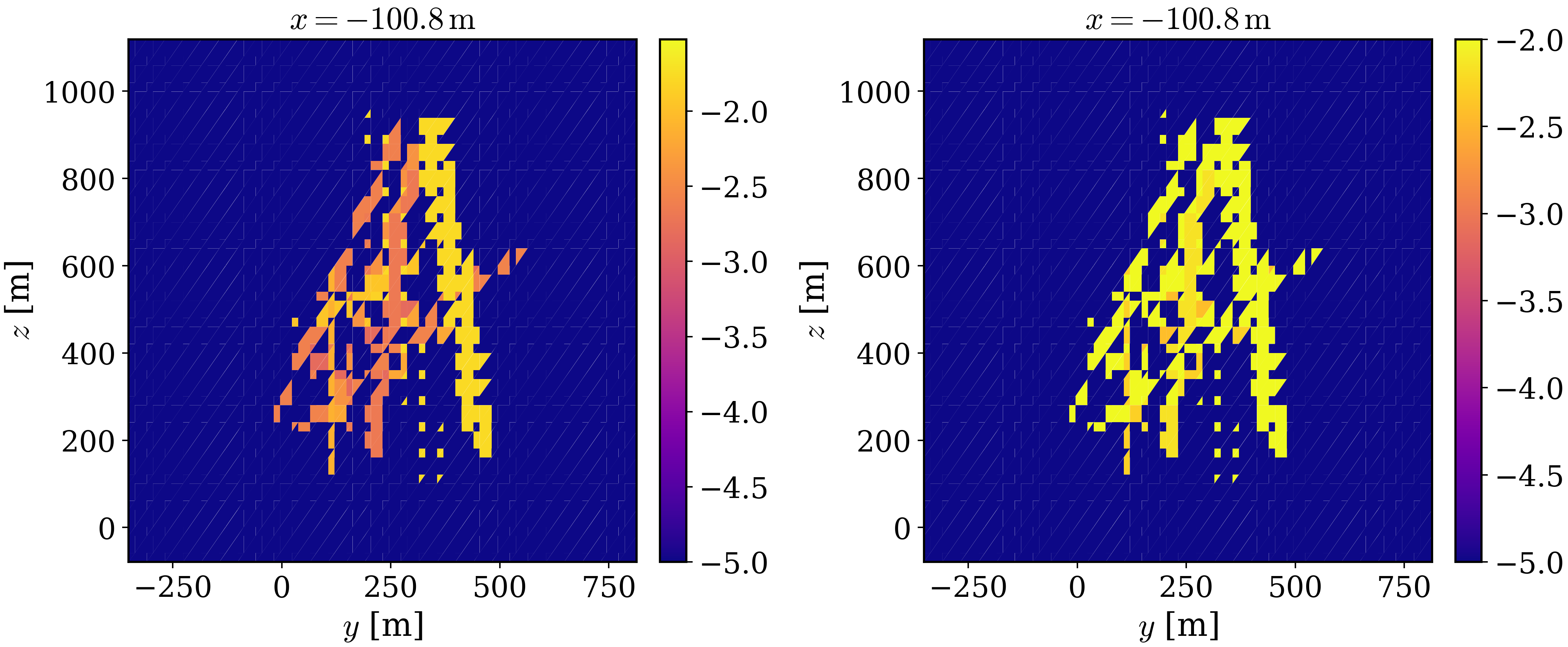}
\caption{Experiment 3. Cross-section at \(x=-100.8\,\mathrm{m}\) of
(a) \(\log_{10}\kappa\) and (b) \(\log_{10}d_T\). The bright regions
represent the fractures, while the dark background represents the
surrounding material with
\(\kappa_{\mathrm{out}}=d_{\mathrm{out}}=10^{-5}\).}
\label{fig:exp3-computational-coefficients-slice}
\end{figure}

We next assign the permeability and thermal conductivity values. Outside fracture supported regions, the permeability and thermal conductivity are fixed at $\kappa_{\mathrm{out}} = d_{\mathrm{out}} = 10^{-5}.$ Each fracture is assigned larger coefficient value according to $\log_{10}\kappa_{f,j} = -2+0.45\xi_j^\kappa,$ and $\log_{10}d_{f,j} = -2+0.35\xi_j^d,$ where $\xi_j^\kappa,\xi_j^d\sim\mathcal{N}(0,1). $ One realization is generated and kept fixed throughout the simulation. If two fracture supported regions overlap, the larger coefficent value is used.  Figures~\ref{fig:exp3-computational-coefficients-3d} and
\ref{fig:exp3-computational-coefficients-slice} show three-dimensional and
two-dimensional views, respectively, of the coefficient fields used in this
experiment.

\subsubsection{Three-dimensional problem setup and fracture-controlled circulation}

We use the same THM parameters, initial conditions, and boundary conditions as in the previous experiments. The problem is solved over $\mathcal{J}=(0,4]$ on a structured $28\times28\times20$ box mesh, with each box divided into tetrahedra. We use continuous piecewise-linear finite elements and the backward Euler method with $\Delta t=0.01$. 

We use four injection locations $\mathbf{x}_1^{\mathrm I},\ldots,\mathbf{x}_4^{\mathrm I}$ and four production locations $\mathbf{x}_1^{\mathrm P},\ldots,\mathbf{x}_4^{\mathrm P}$, as shown in Figure~\ref{fig:exp3-stage8-geometry}. The source rate follows the time-dependent schedule similar to Experiment \ref{subsec:exp2_ghs}, defined as~\eqref{eq:exp3-in-prod-rate}.
Thus, the total injection and production rate remain balanced, while their distributions among the fixed locations vary with time. Because the computational domain is not a unit cube, we scale the spatial distance by its side length. Let $\mathbf{L}=\operatorname{diag}(L_x,L_y,L_z)$, where \(L_x\), \(L_y\), and \(L_z\) are the side lengths of \(\Omega\). Each injection and production location $\mathbf{c}$ is represented by a normalized Gaussian function $\Phi$ centered at $\mathbf{c}$.
$$
\Phi(\mathbf{x};\mathbf{c})
=
\frac{
\exp\left(
-\beta_{\mathrm{src}}
\left\|
\mathbf{L}^{-1}(\mathbf{x}-\mathbf{c})
\right\|_2^2
\right)
}{
\displaystyle
\int_{\Omega}
\exp\left(
-\beta_{\mathrm{src}}
\left\|
\mathbf{L}^{-1}(\mathbf{y}-\mathbf{c})
\right\|_2^2
\right)\,d\mathbf{y}
},
\qquad
\beta_{\mathrm{src}}=1200.
$$
We then set $\Phi_j^{\mathrm I}(\mathbf{x})
=
\Phi(\mathbf{x};\mathbf{x}_j^{\mathrm I})$ and $\Phi_k^{\mathrm P}(\mathbf{x})
=
\Phi(\mathbf{x};\mathbf{x}_k^{\mathrm P}).$ The source terms in the pressure and temperature equations are
$$
\begin{aligned} 
g(\mathbf{x},t)
&=
\sum_{j=1}^{4}
q_j^{\mathrm I}(t)\Phi_j^{\mathrm I}(\mathbf{x})
-
\sum_{k=1}^{4}
q_k^{\mathrm P}(t)\Phi_k^{\mathrm P}(\mathbf{x}),
\\
\eta(\mathbf{x},t)
&=
-\chi_{\mathrm{cold}}
\sum_{j=1}^{4}
q_j^{\mathrm I}(t)\Phi_j^{\mathrm I}(\mathbf{x}),
\end{aligned} $$
where $\chi_{\mathrm{cold}}=0.85$. The source term $g$ represents fluid injection and production, while $\eta$ represents cooling near the injection locations. 

\subsubsection{Geothermal quantities of interest}
We evaluate the same two quantities of interest as interest as in Experiment \ref{subsec:exp2_ghs} with the localized averages now computed over the three dimensional mesh. The total rate $q_{\mathrm{tot}}(t)$ and the allocation weights $w_j^{\mathrm I}(t)$ and $w_k^{\mathrm P}(t)$ are same as those used in the source terms.  Let $\{\mathbf{x}_i\}_{i=1}^{N_v}$ denote the vertices of the three-dimensional mesh. For a scalar finite element function $v_h$ and a location $\mathbf{c}$, we define the local average

$$
\mathcal{A}(v_h;\mathbf{c})
=
\frac{
\displaystyle\sum_{i=1}^{N_v}
\exp\left(
-\beta_{\mathrm{qoi}}
\left\|
\mathbf{L}^{-1}(\mathbf{x}_i-\mathbf{c})
\right\|_2^2
\right)
v_h(\mathbf{x}_i)
}{
\displaystyle\sum_{i=1}^{N_v}
\exp\left(
-\beta_{\mathrm{qoi}}
\left\|
\mathbf{L}^{-1}(\mathbf{x}_i-\mathbf{c})
\right\|_2^2
\right)
},
\qquad
\beta_{\mathrm{qoi}}=1200.$$

Using the localized average, we define the production-region temperature change by 
\begin{equation}
    \overline{\theta}_{\mathrm P}(t^n)
=
\sum_{k=1}^{4}
w_k^{\mathrm P}(t^n)
\mathcal{A}
\left(
\theta_h(\cdot,t^n);
\mathbf{x}_k^{\mathrm P}
\right).
\end{equation}
Similarly, the injection- and production- region pressure are 
$$
    \overline{p}_{\mathrm I}(t^n)
=
\sum_{j=1}^{4}
w_j^{\mathrm I}(t^n)
\mathcal{A}
\left(
p_h(\cdot,t^n);
\mathbf{x}_j^{\mathrm I}
\right),
$$
and
$$
\overline{p}_{\mathrm P}(t^n)
=
\sum_{k=1}^{4}
w_k^{\mathrm P}(t^n)
\mathcal{A}
\left(
p_h(\cdot,t^n);
\mathbf{x}_k^{\mathrm P}
\right).
$$

The two quantities of interest are then defined by
$$
\mathcal{Q}_{\theta}(t^n)
=
-\overline{\theta}_{\mathrm P}(t^n),
\qquad
\mathcal{Q}_{p}(t^n)
=
\frac{
\overline{p}_{\mathrm I}(t^n)
-
\overline{p}_{\mathrm P}(t^n)
}{
q_{\mathrm{tot}}(t^n)
},
$$
where $\mathcal{Q}_{\theta}$ measures the temperature decrease near the production locations, while \(\mathcal{Q}_{p}\) measures the injection-production pressure difference per unit total source rate.

The mesh contains $N_v=17,661$ vertices. Since the displacement has three components, each full THM snapshot has \(N_h=5N_v=88\,305\) entries. We collect \(N_t=400\) snapshots, giving $ \mathbf{Z}\in\mathbb{R}^{88\,305\times400}.$
\subsubsection{Reconstruction results}
We apply SMN using $s=100$ candidate POD modes and $r=10$ selected modes. The selected mode indices are $I_{\mathrm{SMN}} = \{0,3,1,2,4,6,8,5,7,10\}.$ The full-state and component-wise reconstruction errors are reported in Table~\ref{tab:exp3-reconstruction-errors}.

\begin{table}[!h]
\centering
\small
\renewcommand{\arraystretch}{1.10}
\caption{Experiment~3. full-state and component-wise reconstruction errors using \(s=100\) and \(r=10\).}
\label{tab:exp3-reconstruction-errors}
\begin{tabular}{|c|c|c|}
\hline
State or component & POD & SMN \\
\hline
Full THM state
& \(E^{\mathrm{POD}}=4.9372\times10^{-4}\)
& \(E^{\mathrm{SMN}}=1.0712\times10^{-11}\) \\
\(u_x\)
& \(E_{u_x}^{\mathrm{POD}}=2.1562\times10^{-4}\)
& \(E_{u_x}^{\mathrm{SMN}}=3.4983\times10^{-12}\) \\
\(u_y\)
& \(E_{u_y}^{\mathrm{POD}}=1.9807\times10^{-4}\)
& \(E_{u_y}^{\mathrm{SMN}}=3.7173\times10^{-12}\) \\
\(u_z\)
& \(E_{u_z}^{\mathrm{POD}}=2.1237\times10^{-4}\)
& \(E_{u_z}^{\mathrm{SMN}}=3.3616\times10^{-12}\) \\
\(p\)
& \(E_p^{\mathrm{POD}}=2.9676\times10^{-4}\)
& \(E_p^{\mathrm{SMN}}=4.9355\times10^{-12}\) \\
\(\theta\)
& \(E_\theta^{\mathrm{POD}}=5.9076\times10^{-4}\)
& \(E_\theta^{\mathrm{SMN}}=1.3258\times10^{-11}\) \\
\(|\mathbf{u}|\)
& \(E_{|\mathbf{u}|}^{\mathrm{POD}}=1.5249\times10^{-4}\)
& \(E_{|\mathbf{u}|}^{\mathrm{SMN}}=2.0974\times10^{-12}\) \\
\hline
\end{tabular}
\end{table}

Table~\ref{tab:exp3-reconstruction-errors} shows that SMN gives smaller errors than POD for the full THM state and for every solution component. The full-state error is reduced from \(4.9372\times10^{-4}\) to \(1.0712\times10^{-11}\). The SMN component errors are between \(10^{-12}\) and \(10^{-11}\), showing that the selected modes and nonlinear reconstruction accurately represent the three-dimensional snapshot data.

\begin{figure}[!h]
    \centering

    \begin{subfigure}[t]{0.48\textwidth}
        \centering
        \includegraphics[width=\textwidth]
        {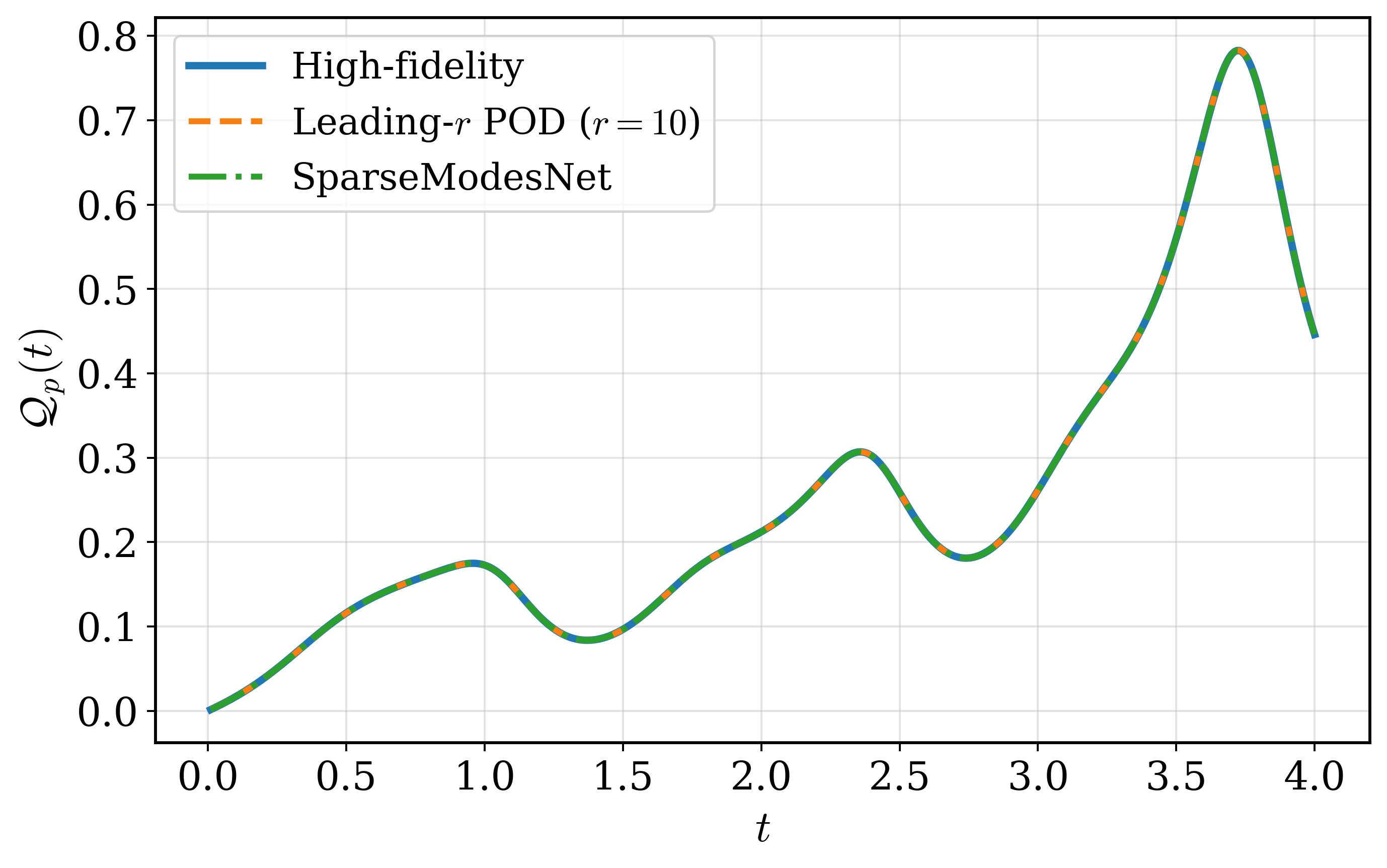}
        \caption{Source-normalized pressure difference 
        $\mathcal{Q}_{p}(t^n)$.}
        \label{fig:exp3-qoi-pressure}
    \end{subfigure}
    \hfill
    \begin{subfigure}[t]{0.48\textwidth}
        \centering
        \includegraphics[width=\textwidth]
        {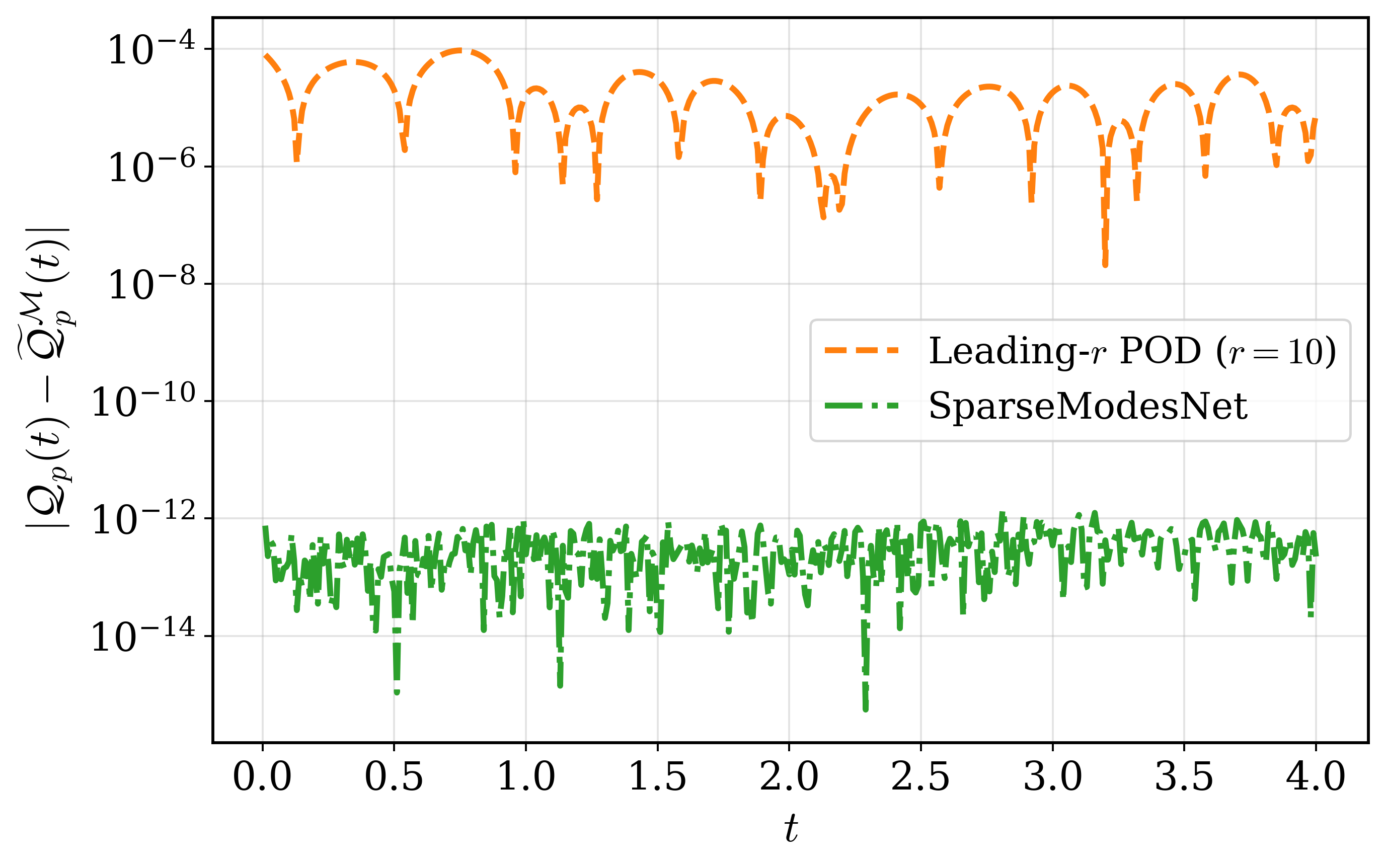}
        \caption{Pointwise absolute errors in 
        $\mathcal{Q}_{p}(t^n)$.}
        \label{fig:exp3-qoi-pressure-error}
    \end{subfigure}

    \caption{Experiment 3. Source-normalized pressure difference and its pointwise
    absolute errors for the leading-$r$ POD and SMN reconstructions,
    using $s=100$ and $r=10$.}
    \label{fig:exp3-qoi-pressure-results}
\end{figure}

\begin{figure}[!h]
    \centering

    \begin{subfigure}[t]{0.48\textwidth}
        \centering
        \includegraphics[width=\textwidth]
        {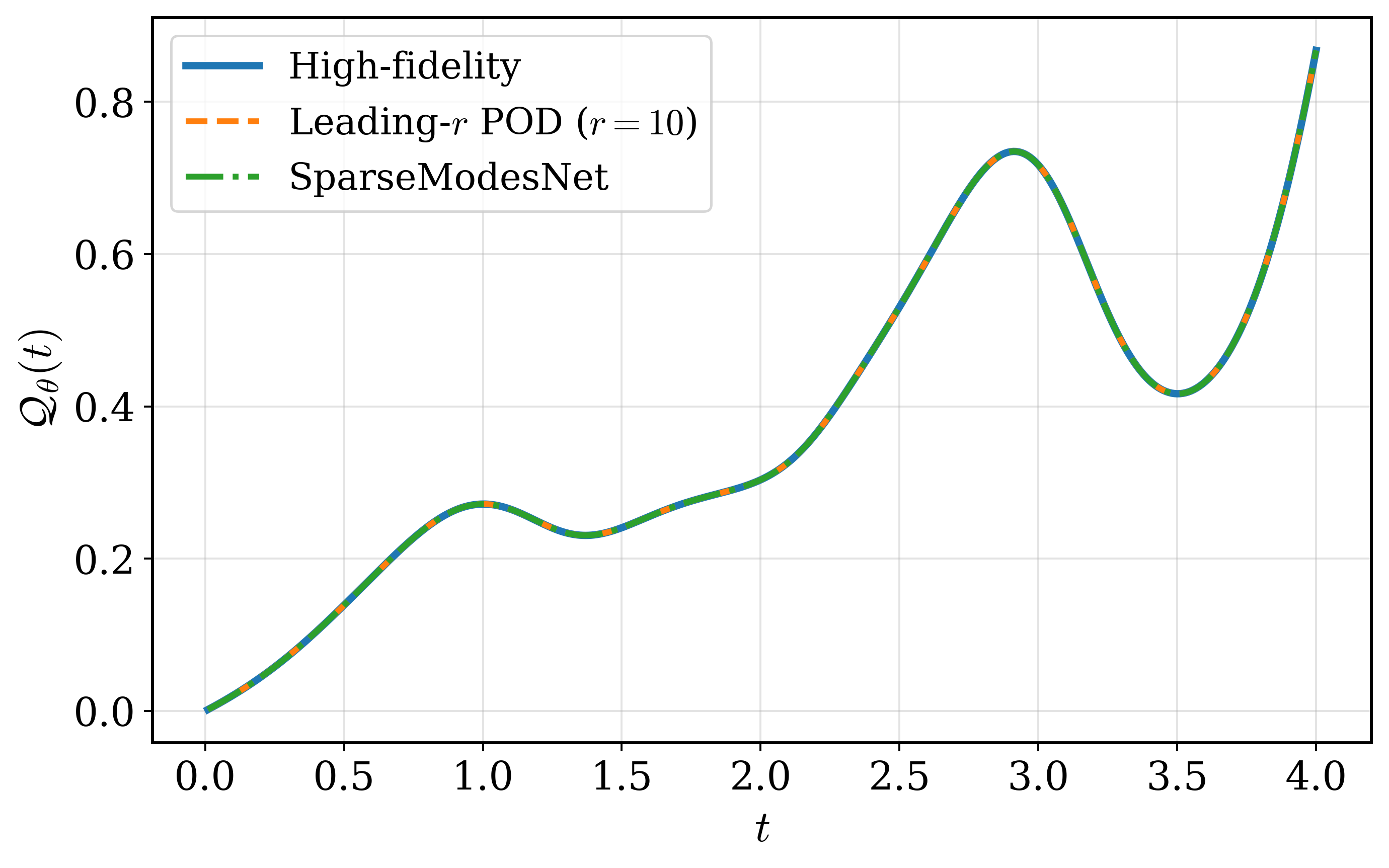}
        \caption{Average temperature near the production
        locations, $\mathcal{Q}_{\theta}(t^n)$.}
        \label{fig:exp3-qoi-temperature}
    \end{subfigure}
    \hfill
    \begin{subfigure}[t]{0.48\textwidth}
        \centering
        \includegraphics[width=\textwidth]
        {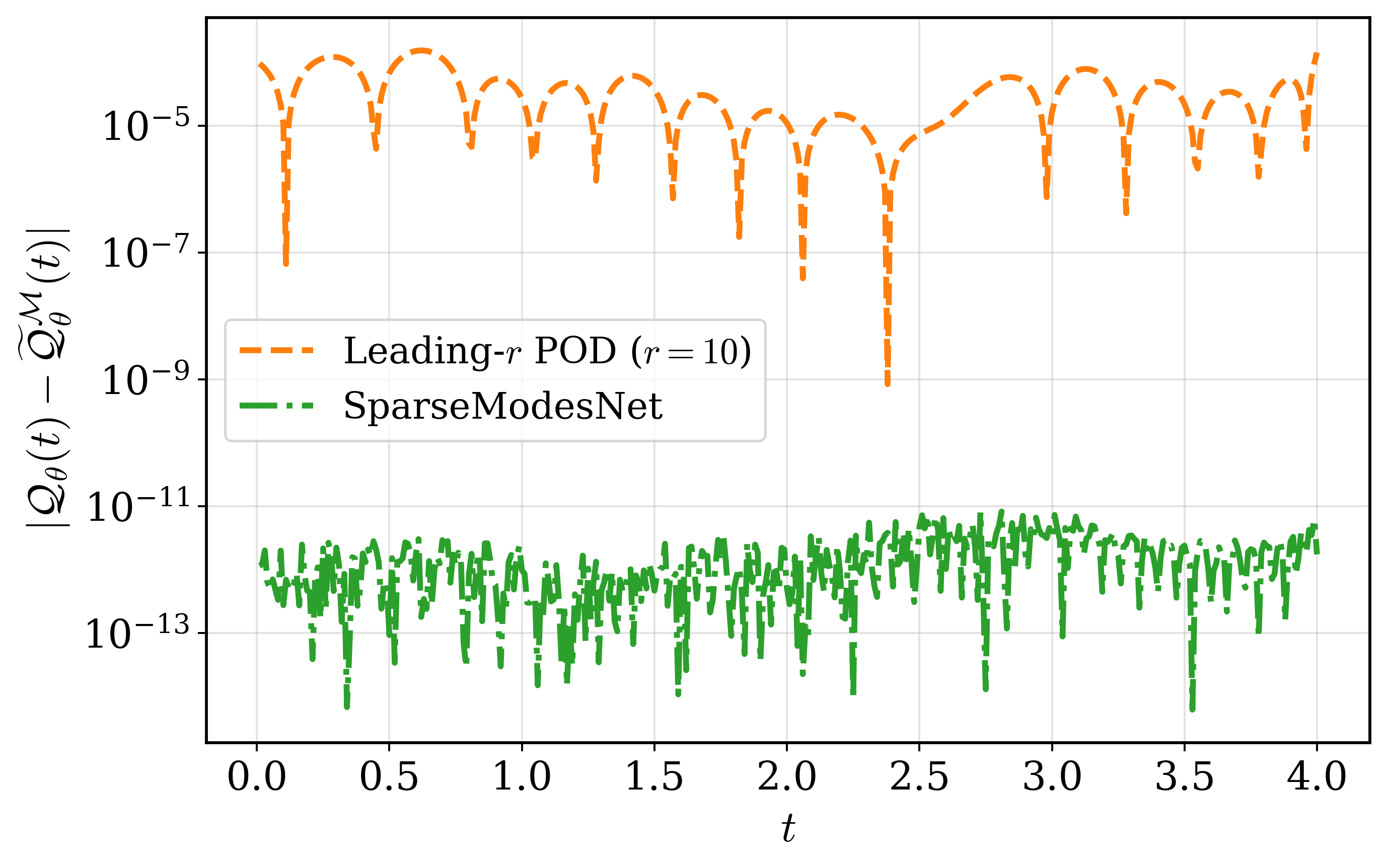}
        \caption{Pointwise absolute errors in 
        $\mathcal{Q}_{\theta}(t^n)$.}
        \label{fig:exp3-qoi-temperature-error}
    \end{subfigure}

    \caption{Experiment 3. Average temperature near the production locations
    and its pointwise absolute errors for the leading-$r$ POD and SMN
    reconstructions, using $s=100$ and $r=10$.}
    \label{fig:exp3-qoi-temperature-results}
\end{figure}

We next examine the two quantities of interest. 
Figures~\ref{fig:exp3-qoi-pressure-results} 
and~\ref{fig:exp3-qoi-temperature-results} compare the quantities of 
interest obtained from the high-fidelity, POD, and SMN fields, together 
with their pointwise absolute errors.

For both quantities of interest, the high-fidelity, POD, and SMN curves
are nearly indistinguishable on the plotted scale. The pointwise
absolute-error plots in Figures~\ref{fig:exp3-qoi-pressure-error}
and~\ref{fig:exp3-qoi-temperature-error}, however, clearly demonstrate
the greater accuracy of the SMN reconstruction. Using the relative error
defined in Experiment~\ref{subsec:exp2_ghs}, the pressure errors are
$E_{\mathcal{Q}_{p}}^{\mathrm{POD}}
    = 1.354 \times 10^{-4},~\text{and}~
E_{\mathcal{Q}_{p}}^{\mathrm{SMN}}
    = 4.324 \times 10^{-12}$.
The corresponding temperature errors are $E_{\mathcal{Q}_{\theta}}^{\mathrm{POD}}
    = 5.489 \times 10^{-4},~\text{and}~
E_{\mathcal{Q}_{\theta}}^{\mathrm{SMN}}
    = 3.645 \times 10^{-11}$.
These results show that both methods preserve the two quantities of
interest, with SMN providing substantially more accurate reconstructions.

\begin{figure}[!h]
\centering
\begin{subfigure}[t]{0.98\textwidth}
\centering
\includegraphics[width=\linewidth]
{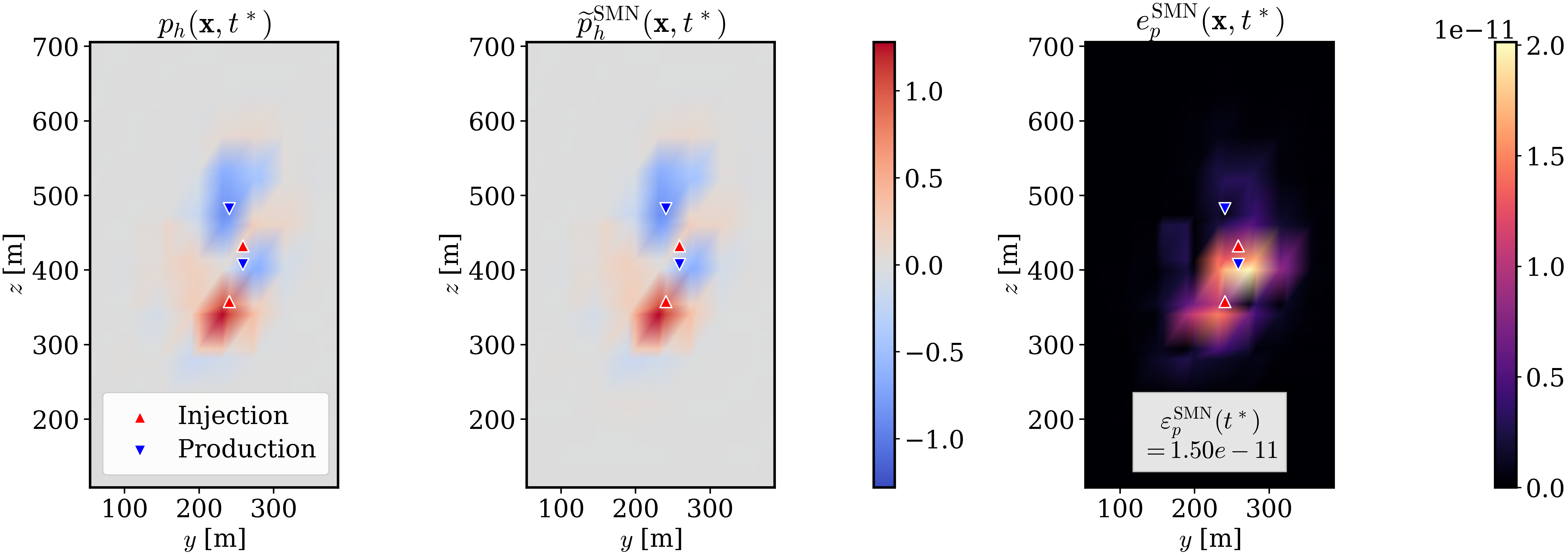}
\caption{Pressure reconstruction.}
\label{fig:exp3-pressure-slice}
\end{subfigure}

\begin{subfigure}[t]{0.98\textwidth}
\centering
\includegraphics[width=\linewidth]
{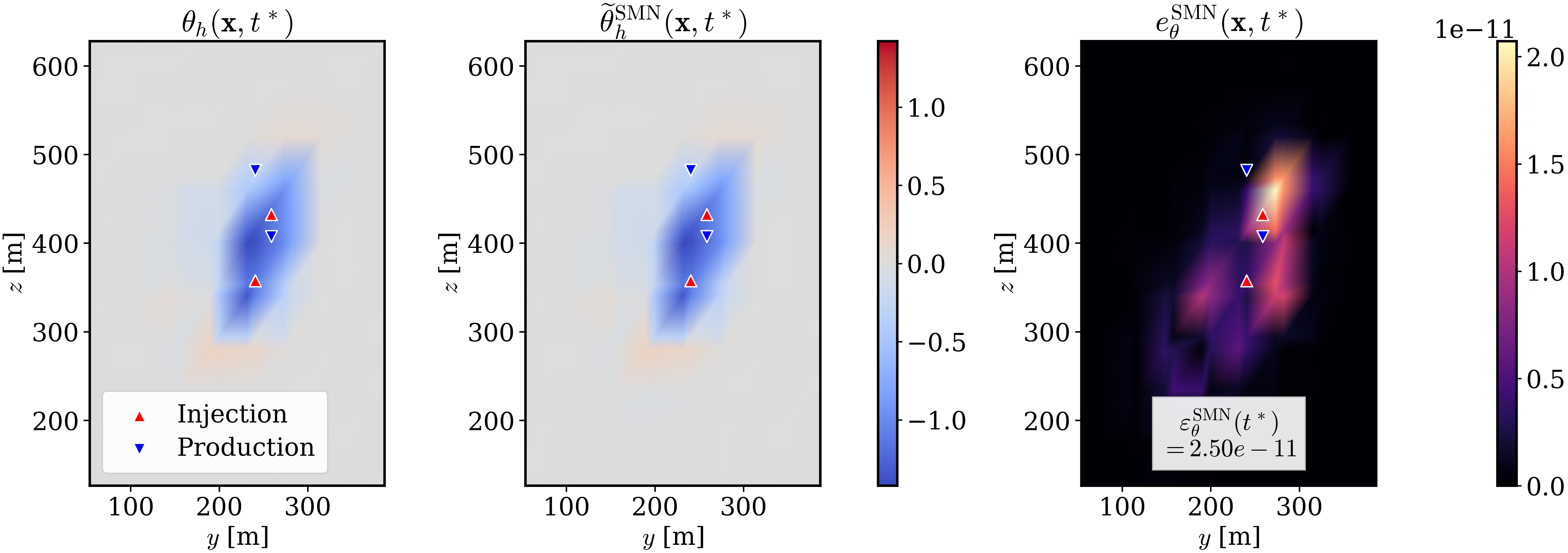}
\caption{Temperature reconstruction.}
\label{fig:exp3-temperature-slice}
\end{subfigure}

\caption{Experiment 3. High-fidelity solutions and SMN reconstructions on the
cross-section \(x=61\,\mathrm{m}\) of the tetrahedral mesh at \(t=2\),
using \(r=10\) selected modes. For each field, the first two panels show the
high-fidelity and reconstructed solutions, while the third panel shows the
pointwise absolute error. The relative errors reported in the text are
computed over the complete three-dimensional nodal fields.}
\label{fig:exp3-pressure-temperature-reconstruction}
\end{figure}

Figure~\ref{fig:exp3-pressure-temperature-reconstruction} compares the high-fidelity pressure and temperature fields with their SMN reconstruction on cross-section $x=61\,\mathrm{m}$ at $t=2$. The relative errors over the complete three-dimensional fields are $\varepsilon_p^{\mathrm{SMN}}=6.67\times10^{-12}$ and $\varepsilon_\theta^{\mathrm{SMN}}=1.44\times10^{-11}.$ SMN reconstructions capture the spatial patterns and heterogeneity between injection and production that are observed on original fields, highlighting the flexibility of SMN.   

\section{Conclusion}
\label{sec:conclusion}
This work presents the first application of SparseModesNet to complex fields arising from coupled multi-physics problems, an area previously unexplored with SparseModesNet.  
Results show that SparseModesNet reconstructs the displacement, pressure,
and temperature fields more accurately than standard POD when both methods use
the same number of reduced coordinates. This improvement is observed across
three experiments of increasing complexity: the two-dimensional heterogeneous problem, the
two-dimensional fracture-network problem, and the three-dimensional Utah FORGE
problem. SparseModesNet also closely preserves the selected pressure- and
temperature-based geothermal quantities of interest. Overall, SparseModesNet
provides an accurate reduced representation of complex THM solution
trajectories in heterogeneous and fractured geothermal systems. 
Such reductions establish a foundation for higher-fidelity reduced-order models, which are essential in complex physics, such as THM problems. Additionally, SparseModesNet demonstrates strong performance in time interpolation, suggesting potential for further exploration in extrapolation and out-of-sample reduction tasks. Preliminary results in these areas indicate that additional parameterization of SparseModesNet may be necessary, though this would expand the method's applicability.

\section*{Acknowledgments}
The authors would like to thank Elizabeth Qian for reviewing the manuscript and providing helpful comments. The work of K. Chawla and S. Lee was supported by the U.S. Department of Energy, Office of Science, Energy Earthshots Initiatives under Award Number DE-SC 0024703.
This work was in part authored by the National Laboratory of the
Rockies for the U.S. Department of Energy (DOE) under
Contract No. DE-AC36-08GO28308. Funding was provided by by the U.S. Department of Energy, Office of Science, Energy Earthshots Initiatives under Award Number DE-SC 0024703. The views expressed do not necessarily
represent the views of the U.S. Department of Energy or
the U.S. Government. The U.S. Government retains and
the publisher, by accepting the article for publication,
acknowledges that the U.S. Government retains a
nonexclusive, paid-up, irrevocable, worldwide license to
publish or reproduce the published form of this work, or
allow others to do so, for U.S. Government purposes. NLR/JA-2C00-101929.

\bibliographystyle{elsarticle-num}
\bibliography{ref}

\end{document}